\documentclass[12pt]{article}

\usepackage[margin=1in]{geometry}  
\usepackage{graphicx}              
\usepackage{amsmath}               
\usepackage{amsfonts}              
\usepackage{amsthm}                
\usepackage{amssymb}
\usepackage{amscd}
\usepackage{amsmath,amscd}
\usepackage{float}
\usepackage{url}

\usepackage{graphics}
\usepackage{pstricks, pst-node, pst-tree, pst-math, pst-xkey}

\usepackage[labelformat=empty]{caption}

\begin{document}

\nocite{*}

\title{ACC for minimal lengths of extremal rays for surfaces.}

\author{Evgeny Mayanskiy}

\maketitle

\begin{abstract}
  We prove that the lengths of extremal rays of log canonical Fano {\it surfaces} with Picard number one satisfy the ascending chain condition. This confirms the $2$-dimensional case of a conjecture stated by Fujino and Ishitsuka in \cite{Fujino}.
\end{abstract}

\section{Introduction.}

Fix $n\geq 2$ and a DCC set $I\subset \mathbb R$. Consider the set ${\mathcal F}_n(I)$ of all log canonical pairs $(X,\Delta)$, where $X$ is a normal $\mathbb Q$-factorial projective variety of dimension $n$ and Picard number one, the coefficients of the $\mathbb Q$-divisor $\Delta$ are in $I$ and $-(K_X+\Delta)$ is ample.\\

Given $(X,\Delta)\in {\mathcal F}_n(I)$, define $l(X,\Delta)>0$ to be the minimal intersection number $-(K_X+\Delta)\cdot C$ for curves $C$ on $X$.\\

Fujino and Ishitsuka formulated the following conjecture in \cite{Fujino}.\\

{\bf Conjecture (cf. \cite{Fujino}).} {\it Given $n\geq 2$ and a DCC set $I$, the set
$$
\{ \; l(X,\Delta) \;\mid \; (X,\Delta)\in {\mathcal F}_n(I)\; \} \subset \mathbb R
$$
satisfies ACC.}\\
 
This Conjecture was verified for toric Fano varieties with $I=\varnothing$ in \cite{Fujino}. The purpose of this note is to check this Conjecture in the case $n=2$ and $I=\varnothing$.\\

Let $X$ be a normal projective surface with log canonical singularities such that $-K_X$ is ample and the Picard number $\rho(X)$ is $1$. Denote by $l(X)>0$ the minimal intersection number $-K_X\cdot C$ for curves $C$ on $X$. Let ${\mathcal F}_2$ be the set of all such surfaces $X$.\\

The main result of this note is the following Theorem.\\

{\bf Theorem 2.} {\it The set $\{\; l(X)\; \mid \; X \in {\mathcal F}_2 \; \} \subset \mathbb R $ satisfies ACC.}\\

Our proof of Theorem~2 follows the same scheme as the proof by Alexeev of a similar result on fractional indices of log del Pezzo surfaces \cite{Alexeev} and depends significantly on some ideas of Nikulin \cite{NikulinII}, \cite{NikulinI}.\\

We use freely the notation and terminology from \cite{KM}. ACC means 'ascending chain condition', DCC means 'descending chain condition'. Given a normal projective surface $X$, we denote by $f\colon Y\to X$ its minimal resolution and by $\pi\colon Y\to Z$ a minimal model. $\rho(X)$ denotes the Picard number of $X$.\\

In Section~2 we verify certain bounds on the Picard numbers of minimal resolutions of surfaces with ample anti-canonical bundles. These bounds were established earlier by Alexeev \cite{Alexeev}. In order to verify them, we use a natural generalization of log terminal systems studied by Nikulin (see \cite{NikulinII}). This generalization is spelled out in the Appendix. In Section~4 we apply the Picard number bounds to prove our main Theorem~2. Section~3 is devoted to examples. \\

\section{Bounds on the Picard numbers of del Pezzo surfaces.}

In this Section we study bounds on the Picard numbers of certain surfaces largely following the paper \cite{Alexeev} by Alexeev. We do not prove any new results here.\\

{\bf Definition.} For any $\epsilon >0$ let $DP_{\epsilon}$ denote the set of all normal projective log canonical surfaces $X$ such that $-K_X$ is ample and moreover $-K_X\cdot C \geq \epsilon$ for every curve $C$ on $X$. (This is a variant of the Condition $DP(\epsilon)$ of Alexeev from \cite{Alexeev},~2.1.)\\

{\bf Remark.} If $X$ is a normal \textit{non-rational} projective log canonical surface with ample $-K_X$, then $X$ is the contraction of the section with negative self-intersection of $Y={\mathbb P}_E(\mathcal E)\to E$, where $\mathcal E$ is a normalized rank $2$ locally free sheaf on a smooth elliptic curve $E$ and the parameter $e=-deg(\mathcal E)\geq 2$. In particular, $\rho(Y)=2$ for the minimal resolution $Y$ of such a surface $X$. (See \cite{Badescu}, Theorem~2 or \cite{Fujisawa}, Theorem~2.1.)\\

For any $X\in DP_{\epsilon}$ (with the minimal resolution $f\colon Y\to X$) one constructs the ampleness polytope $\mathcal M(Y)$ of $Y$ in a hyperbolic space $\mathcal L(Y)$ as follows (see \cite{NikulinI}).\\

Consider the cone $V(Y)=\{ \; x\in NS(Y)\otimes_{\mathbb Z} \mathbb R \; \mid \; x^2>0 \; \}$ in the vector space generated by the Neron-Severi group $NS(Y)$ of $Y$. Take $V^{+}(Y)$ to be the connected component of $V(Y)$ containing an ample divisor class and define
$$
\mathcal L(Y)=V^{+}(Y)/{\mathbb R}^{+},
$$
where ${\mathbb R}^{+}$ is the group of positive reals acting on $V^{+}(Y)$ by multiplication.\\

For any $\delta\in NS(Y)\otimes_{\mathbb Z} \mathbb R$ with ${\delta}^2<0$ one defines the half-space
$$
{\mathcal H}^{+}_{\delta}=\{ \; {\mathbb R}^{+}\cdot x\in \overline{\mathcal L (Y)} \;\mid \; x\cdot \delta \geq 0 \; \}.
$$

The ampleness polytope of $Y$ is defined as 
$$
\mathcal M (Y)=\bigcap {\mathcal H}^{+}_{\delta} \subset \overline{\mathcal L (Y)},
$$
where the intersection is taken over the set of all integral curves $\delta$ on $Y$ with negative self-intersection.\\

Assume that $\rho(Y)\geq 3$. Then by \cite{NikulinI}, Lemma~2.2 and Lemma~A.3.0, $\mathcal M (Y)$ is a convex polytope with acute angles, finitely many faces and finite volume.\\

Such polytopes are interesting to us because of the following result.\\

\textbf{Combinatorial Theorem (\cite{NikulinI}, Lemma~3.4)} \textit{Let $\mathcal L$ be a hyperbolic space and $\mathcal M\subset \mathcal L$ a convex polytope with acute angles, finite volume and finitely many faces such that its Gram graph $\Gamma (\mathcal M)$ (formed by the vectors orthogonal to the facets) satisfies the following conditions:
\begin{enumerate}
\item[a)] Every Lanner subgraph of $\Gamma (\mathcal M)$ has at most $l$ vertices;
\item[b)] For every connected elliptic subgraph of $\Gamma (\mathcal M)$ with $n$ vertices, the number of pairs of its vertices with distance $d\leq l-2$ between them is at most $C_1\cdot n$, and the number of pairs of its vertices with distance $d$ between them, where $l-2<d\leq 2l-3$, is at most $C_2\cdot n$.
\end{enumerate}
Then $dim(\mathcal L)\leq 96(C_1+C_2/3)+68$.}\\

The following Theorem was proved in \cite{Alexeev}, Theorem~2.3 (see \cite{Alexeev92}, Remark~5.5). By a remark in the beginning of this Section we may assume in the proof that the surface $X\in DP_{\epsilon}$ is rational (and more generally that $\rho(Y)\geq 3$) throughout the proof.\\

{\bf Theorem 1 (Alexeev, cf. \cite{Alexeev}, Theorem 2.3, \cite{Alexeev92}, Remark~5.5).} {\it There in an integer $A>0$ such that for any $\epsilon >0$ and any $X\in DP_{\epsilon}$, the Picard number of the minimal resolution $f\colon Y\to X$ satisfies the bound
$$
\rho(Y)\leq A\cdot (1+1/\epsilon).
$$}\\

{\it Proof:} Alexeev's argument in \cite{Alexeev} goes through in our case. However, in order to apply it literally, we need a slight generalization of Nikulin's theory of log terminal systems (\cite{NikulinII}). The notion we need is that of 'log canonical systems'. Their theory can be developed parallelly to the theory of log terminal systems of Nikulin and is spelled out in the Appendix.\\

The argument in \cite{Alexeev} is based on a bound on the dimension of the hyperbolic space admitting a polytope with certain properties. Namely, it uses the Combinatorial Theorem quoted above, where the hyperbolic space $\mathcal L=\mathcal L(Y)$ and the polytope $\mathcal M=\mathcal M(Y)\cap \mathcal L(Y)\subset \mathcal L$ were defined above.\\ 

Combinatorial Theorem says in this case that certain properties of the Gram graph of the ampleness polytope ${\mathcal M}(Y)$ of $Y$ imply an upper bound on the Picard number $\rho (Y)$ (see \cite{NikulinI} and \cite{Alexeev}, Theorem~2.2).\\

In \cite{Alexeev} Alexeev checked these properties in the case, when $X$ has log terminal singularities (and satisfies the Condition $DP(\epsilon)$), \cite{Alexeev}, Theorem~3.13 and Theorem~3.15.\\

The same argument works in our situation. As we have already remarked, one may assume that $X\in DP_{\epsilon}$ is rational. Hence it has only rational singularities (by \cite{Fujisawa}, Corollary~1.9).\\

One uses the results of the Appendix (a generalization of Nikulin's theory of log terminal systems of vectors from \cite{NikulinII} to the log canonical case) and checks that Alexeev's proof of \cite{Alexeev}, Theorem~3.13 and Theorem~3.15 goes through in our case as well. This implies that the Gram graph of the ampleness polytope ${\mathcal M}(Y)$ of $Y$ enjoys the properties needed to derive the bound on the Picard number $\rho (Y)$ from the Combinatorial Theorem. {\it QED}\\

We will derive Theorem~2 from Theorem~1 and the alternative in the following Corollary, which was proved by Alexeev in \cite{Alexeev}, \cite{AlexeevCorr}.\\

{\bf Corollary (Alexeev, \cite{AlexeevCorr}, Lemma 4.1, \cite{Alexeev}, Lemma 4.2).} {\it Fix $\epsilon >0$. Then for any $N>0$ large enough (for the chosen $\epsilon$) there is $M>0$ with the following property.\\

Given $X\in DP_{\epsilon}$ with the minimal resolution $f\colon Y\to X$, whose exceptional locus consists of the integral curves $F_1$, ..., $F_r$, we have
$$
\rho(Y)\leq N \;\;\; \mbox{and either}
$$
\begin{itemize}
\item[(1)] $-F_i^2\leq M$ for all $i$; or
\item[(2)] there is a birational morphism $Y\to {\mathbb F}_n$ for some $n\geq N$ (here ${\mathbb F}_n={\mathbb P}(\mathcal O \oplus \mathcal O (n))$ denotes a Hirzebruch surface) such that any integral curve on $Y$ with negative self-intersection is either the strict transform of the $(-n)$-curve on ${\mathbb F}_n$, or maps to a fiber of the fibration ${\mathbb F}_n\to {\mathbb P}^1$; or
\item[(3)] $Y={\mathbb P}_E(\mathcal E)$, where $\mathcal E$ is a normalized rank $2$ locally free sheaf on a smooth elliptic curve $E$ and the parameter $e=-deg(\mathcal E)\geq N$, $r=1$ and $F_1$ is the section with negative self-intersection of the fibration ${\mathbb P}_E(\mathcal E)\to E$.
\end{itemize}
}

\section{Examples.}

Before we continue to the proof of Theorem~2, let us consider two examples.\\

{\bf Example 1.} Take $n\geq 2$. Let $X_n\in {\mathcal F}_2$ be the surface obtained by contracting the $(-n)$-curve of the Hirzebruch surface ${\mathbb F}_n$:
$$
f\colon Y={\mathbb F}_n \to X_n.
$$

Denote by $L$ the general fiber of the fibration $g\colon {\mathbb F}_n \to {\mathbb P}^1$ and by $B$ the $(-n)$-section of $g$.\\

Note that 
$$
K_Y=(-2)\cdot B + (-n-2)\cdot L
$$
and
$$
f^{*}(K_{X_n})=K_Y+\left(1-\frac{2}{n}\right)\cdot B.
$$
\par
Then $l(X_n)=-(f^{*}K_{X_n})\cdot L=2-\left(1-\frac{2}{n}\right)=1+\frac{2}{n}$.\\

We see that the set $\{ \; l(X_n)\;\mid \; n\geq 2 \; \}$ satisfies ACC and does not satisfy DCC.\\

{\bf Example 2.} Let $E$ be a smooth elliptic curve, $e\geq 2$ and $\mathcal E$ a normalized rank $2$ locally free sheaf on $E$ such that $deg(\mathcal E)=-e$. Let $X_e\in {\mathcal F}_2$ be the surface obtained by contracting the section with negative self-intersection of $Y={\mathbb P}_E(\mathcal E)\to E$:
$$
f\colon Y={\mathbb P}_E(\mathcal E) \to X_e.
$$

Denote by $L$ the general fiber of the fibration $g\colon {\mathbb P}_E(\mathcal E) \to E$ and by $B$ the $(-e)$-section of $g$.\\

Note that 
$$
f^{*}(K_{X_e}) = K_Y + B.
$$

Then $l(X_e)=-(f^{*}K_{X_e})\cdot L=1$.\\

We see that the set $\{ \; l(X_e)\; \mid \; E,\; \mathcal E,\; e=-deg(\mathcal E)\geq 2  \; \}$ satisfies ACC.\\

\section{Proof of Theorem~2.}

Fix $\epsilon >0$ and take $N>0$ large enough (for Corollary to hold).\\

Let $DP^{1}_{\epsilon}\subset DP_{\epsilon}$ be the subset consisting of \textit{rational} $X\in DP_{\epsilon}$ with Picard number $\rho(X)=1$ and $DP^{f}_{\epsilon}\subset DP^{1}_{\epsilon}$ be the subset consisting of $X\in DP_{\epsilon}$ for which alternative $(2)$ of Corollary holds (and $\rho(X)=1$).\\

Given $X\in DP^{f}_{\epsilon}$, we will denote by $\pi\colon Y\to Z={\mathbb F}_n$ the birational morphism from Corollary, $g\colon Y\to {\mathbb P}^{1}$ the resulting ${\mathbb P}^{1}$-fibration and $L$ the general fiber of $g$.\\

Consider the following sets:
$$
\mathcal A = \{ \;  K_Y^2 \;\mid \; X\in DP_{\epsilon} \; \}, 
$$

$$
\mathcal B_0 =  \{ \;  K_X^2 \;\mid \; X\in DP^{1}_{\epsilon} \setminus DP^{f}_{\epsilon} \; \}, 
$$

\begin{multline*}
\mathcal B =  \{ \;  -(f^{*}K_X)\cdot C \;\mid \; X\in DP^{f}_{\epsilon},\; C\; \mbox{is an integral}\\
\mbox{ curve on}\; Y, C^2<0 \; \mbox{or}\; C \; \mbox{is a fiber of}\; g  \; \}. 
\end{multline*}

{\bf Lemma 3.} {\it (1) $\mathcal A$ and $\mathcal B_0$ are finite;\\

(2) $\mathcal B$ satisfies ACC.}\\

{\it Proof:} (1) Let $X\in DP_{\epsilon}$. Consider the codiscrepancy formula:
$$
f^{*}K_X=K_Y+\sum_i{{\alpha}_i\cdot F_i},
$$
where $F_i$ are the irreducible components of the exceptional locus of the minimal resolution $f\colon Y\to X$. The coefficients ${\alpha}_i\in [ 0,1 ] \cap {\mathbb Q}$ are determined from the linear system:
$$
\sum_i{{\alpha}_i\cdot (F_i\cdot F_j)}=-K_Y\cdot F_j.
$$

It follows from Corollary that as $X\in DP^{1}_{\epsilon}\setminus DP^{f}_{\epsilon}$, the integers $(F_i\cdot F_j)$, $(K_Y\cdot F_i=-F_i^2-2)$ belong to a finite subset of $\mathbb Z$ (determined by $\epsilon$ and $N$).\\

Hence as $X\in DP^{1}_{\epsilon}\setminus DP^{f}_{\epsilon}$, the coefficients ${\alpha}_i$ belong to a finite set of rationals.\\

Since $K_X^2=K_Y\cdot (K_Y+\sum_i{{\alpha}_i\cdot F_i})=K_Y^2+\sum_i{{\alpha}_i\cdot (K_Y\cdot F_i)}$, the finiteness of $\mathcal A$ implies the finiteness of $\mathcal B_0$.\\

Let us show that $\mathcal A$ is finite.\\

We may assume that $X\in DP_{\epsilon}$ is rational. By Noether's formula, $K_Y^2=10-\rho(Y)$ lies in a finite set of integers (determined by $\epsilon$). Hence $\mathcal A$ is finite.\\ 

(2) Let $X\in DP^{f}_{\epsilon}$. Let $F_0$ be the strict transform (via $\pi$) on $Y$ of the $(-n)$-curve on $Z={\mathbb F}_n$ and $m=-F_0^2 > 0$.\\

Let $F_1$, ..., $F_p$ be the components of the fibers of $g$ with self-intersection $F_i^2\leq -2$. Then $f$ contracts exactly $F_0$, $F_1$, ..., $F_p$.\\

Consider the codiscrepancy formula:
$$
f^{*}K_X=K_Y+\sum_{i=0}^{p}{{\alpha}_i\cdot F_i},
$$
where ${\alpha}_i\in [ 0,1 ] \cap {\mathbb Q}$.\\

The coefficients ${\alpha}_i$ can be computed from the following system of linear equations:
$$
A\cdot \left(
\begin{array}{c}  
{\alpha}_0 \\
{\alpha}_1 \\
\ldots \\
{\alpha}_p
\end{array}
\right) = \left(
\begin{array}{c}  
{\beta}_0 \\
{\beta}_1 \\
\ldots \\
{\beta}_p
\end{array}
\right) ,
$$ 
where $A=\left( -(F_i\cdot F_j) \right)_{i,j}$ is (the negative of) the intersection matrix of the curves $F_i$ and ${\beta}_j=K_Y\cdot F_j=-F_j^2-2$ for every $j$.\\

We know from Corollary that (assuming $\epsilon$ and $N$ fixed) the set of all possible values of $p$, ${\beta}_1$, ..., ${\beta}_p$, $F_i\cdot F_j$, $F_0\cdot F_j$, $i,j\geq 1$, is a \textit{finite} set of integers.\\

By substituting $F_1$, ..., $F_p$ with their suitable ${\mathbb R}$-linear combinations, we may rewrite $A$ in the following form:
$$
A = \left(
\begin{array}{ccccc}  
m & a_1 & \ldots & \ldots & a_p \\
a_1 & 1 &        & & \\
\vdots &  & \ddots & 0 & \\
\vdots & 0 &  & \ddots & \\
a_p &  &   & & 1 
\end{array}
\right)
$$
(i.e. we diagonalized the symmetric positive definite matrix $\left( -F_i\cdot F_j\right)_{i,j\geq 1}$).\\

Hence the coefficients $({\alpha}_0, {\alpha}'_1, ...{\alpha}'_p)$ in the respective codiscrepancy formula form the solution of the linear system
$$
\left(
\begin{array}{ccccc}  
m & a_1 & \ldots & \ldots & a_p \\
a_1 & 1 &        & & \\
\vdots &  & \ddots & 0 & \\
\vdots & 0 &  & \ddots & \\
a_p &  &   & & 1 
\end{array}
\right) \cdot \left(
\begin{array}{c}  
{\alpha}_0 \\
{\alpha}'_1 \\
\ldots \\
{\alpha}'_p
\end{array}
\right)  = \left(
\begin{array}{c}  
m-2 \\
b_1 \\
\ldots \\
b_p
\end{array}
\right) .
$$

Here $a_i$, $b_j$ run over a finite set of reals (as $X$ runs over $DP^{f}_{\epsilon}$).\\

Since
$$
A^{-1}=\frac{1}{m-a_1^2- ... -a_p^2}\cdot \left(
\begin{array}{cccc}  
1 & -a_1 & -a_2 & \vdots \\
-a_1 & m-a_2^2- ... -a_p^2 & a_1a_2 & \vdots \\
-a_2 & a_1a_2 & m-a_1^2-a_3^2- ... -a_p^2 & \vdots \\
-a_3 & a_1a_3 & a_2a_3 & \vdots \\
\ldots & \ldots & \ldots & \vdots \\
-a_p & a_1a_p & a_2a_p & \vdots 
\end{array}
\right) ,
$$

we compute:
$$
{\alpha}_0 =\frac{m-2-a_1b_1-...-a_pb_p}{m-a_1^2- ... -a_p^2} \in [ 0,1 ] .
$$

By Corollary every integral curve $C$ on $Y$ with negative self-intersection is either one of $F_i$ or is a $(-1)$-curve in a fiber of $g$. In the former case $-(f^{*}K_X)\cdot C=0$.\\

Since $\rho(X)=1$, every reducible fiber of $g$ contains exactly one $(-1)$-curve. Upto adding a linear combination of $F_1$, ..., $F_p$, such a curve $C$ is numerically equivalent to $\lambda \cdot L$, where $\lambda \in \mathbb Q$ is an element of a finite set of rationals (determined by $\epsilon$ and $N$).\\

In particular, $-(f^{*}K_X)\cdot C=\lambda \cdot (-f^{*}K_X)\cdot L=\lambda \cdot (2-{\alpha}_0)$.\\

From the expression for ${\alpha}_0$ above we see that ${\alpha}_0\to 1$ as $m\to \infty$. Hence the set of all possible values of ${\alpha}_0\in [ 0,1 ]$ satisfies DCC.\\

We conclude that ${\mathcal B}$ satisfies ACC. {\it QED}\\

Now we are ready to prove the main result of this note.\\

{\bf Theorem 2.} {\it Conjecture holds for $n=2$ and $I=\varnothing$.}\\

{\it Proof:} Arguing by contradiction, it is enough to show that for any $\epsilon > 0$, the set $\{  \; l(X) \; \mid \; X\in DP^{1}_{\epsilon} \; \} $ satisfies ACC.\\

It is clear that the set $\{  \; l(X) \; \mid \; X\in DP^{1}_{\epsilon} \setminus DP^{f}_{\epsilon} \; \} $ satisfies ACC. Indeed, for any $\epsilon >0$ there is $P\in \mathbb Z$ such that $P\cdot f^{*}(K_X)$ is Cartier for any $X\in DP^{1}_{\epsilon} \setminus DP^{f}_{\epsilon}$. In particular, $l(X)\in \frac{1}{P}\cdot \mathbb Z$ for any $X\in DP^{1}_{\epsilon} \setminus DP^{f}_{\epsilon}$. On the other hand, $0< l(X) \leq 4$ (see \cite{Fujino2}, Theorem~18.2). Hence the set $\{  \; l(X) \; \mid \; X\in DP^{1}_{\epsilon} \setminus DP^{f}_{\epsilon} \; \} $ is finite.\\

Let $X\in DP^{f}_{\epsilon}$ and $D$ be the reduced exceptional divisor of the minimal resolution $f\colon Y\to X$.\\

If $l$ and $b$ are the general fiber and the $(-n)$-section of the fibration $Z={\mathbb F}_n\to {\mathbb P}^1$, then
$$
K_Z=(-2)\cdot b+(-n-2)\cdot l.
$$

Let $C$ be an integral curve on $Y$ such that $(-f^{*}K_X)\cdot C = l(X)$.\\

{\bf Claim.} If $N$ in Corollary is taken large enough (for a given $\epsilon$), then $\mid K_Y+D+C \mid = \varnothing $.\\

{\it Proof of the Claim:} $K_Y={\pi}^{*}K_Z+\Gamma$ for some $\Gamma \geq 0$ supported on reducible fibers of $g$. Hence 
$$
K_Y\leq -2\cdot {\pi}^{*}(b)-N\cdot {\pi}^{*}(l)+\Gamma .
$$

Since the coefficients of the irreducible components of $\Gamma$ belong to a finite set of integers (depending on $\epsilon$), one can take $N$ in Corollary large enough (for the given $\epsilon$) so that 
$$
K_Y\leq -2\cdot {\pi}^{*}(b)-\left( 2+\rho(Y)+\frac{N}{2} \right) \cdot L.
$$
If $K_Y+D+C \geq 0$, then 
$$
K_Y+B+\rho(Y)\cdot L+C \geq 0,
$$
where $B$ is the strict transform of $b$ on $Y$.\\

Hence $B+C-2\cdot {\pi}^{*}(b)-\left( 2+\frac{N}{2} \right) \cdot L \geq 0$. This implies that $C\geq 2L$.\\

Hence $l(X)=(-f^{*}K_X)\cdot C \geq 2\cdot (-f^{*}K_X)\cdot L > (-f^{*}K_X)\cdot L \geq l(X)$, which is absurd.\\

{\it The Claim is proved.}\\

So, we may assume that $\mid K_Y+D+C \mid = \varnothing $.\\

By \cite{KojimaTakahashi}, Lemma~3.6 or \cite{Zhang88}, Lemma~2.2 either $\rho(Y)=2$ (i.e. $X$ is obtained by contracting the $(-n)$-curve on $Y={\mathbb F}_n$) or $C$ can be taken to be a smooth rational curve with self-intersection $C^2=-1$.\\

By Lemma~3(2) and Example~1 we conclude that the set
$$
\{  \; l(X) \; \mid \; X\in DP^{f}_{\epsilon} \; \}
$$
satisfies ACC. Together with Example~2 this implies the Theorem. {\it QED}\\

\newpage

\section*{Appendix. Log canonical systems of vectors (following \cite{NikulinII})}

In \cite{NikulinII} Nikulin presents the theory of log terminal systems of vectors which, in particular, can be applied to describe configurations of curves with negative self-intersection on minimal resolutions of surfaces (of a certain kind) with log terminal singularities (see, for example, \cite{NikulinII}, \cite{Alexeev}). In fact, this theory is also applicable to surfaces with log canonical singularities. The purpose of this Appendix is to spell out such an extension.\\

We follow closely the exposition of \cite{NikulinII} and do not claim any originality. We left out proofs whenever Nikulin's arguments in \cite{NikulinII} go through.\\

{\bf Definition (\cite{NikulinII}).} A finite set of vectors $\mathcal V =\{ v_1, ..., v_n \}$ of a finite-dimensional quadratic module over $\mathbb Q$ (with bilinear form $(u,v)\mapsto u\cdot v$) is called an {\it at most hyperbolic system of vectors}, if 
\begin{itemize}
\item its {\it Gram matrix} $A=(a_{ij})$, $a_{ij}=v_i\cdot v_j$, satisfies: $a_{ij}\in \mathbb Z$, $a_{ii}=-b_i\leq -1$ for any $i,j$ and $a_{ij}\geq 0$ if $i\neq j$; and
\item the lattice $S$ determined by $A$ (by modding out the kernel of $A$) has signature $(p,q)$ with $p\leq 1$.
\end{itemize}

{\bf Remarks (\cite{NikulinII}).} (1) At most hyperbolic systems of vectors can be encoded by weighted graphs (with the Gram matrix $A$ of an at most hyperbolic system of vectors $\mathcal V$ being the Gram matrix of the {\it associated weighted graph}).\\

(2) One can speak naturally of {\it equivalent} at most hyperbolic systems of vectors and a {\it subsystem} of an at most hyperbolic system of vectors.\\

(3) We will abuse the language by calling an 'at most hyperbolic system of vectors' simply a 'system'.\\ 

(4) If $A=(a_{ij})$ is a negative definite symmetric square matrix with real entries such that $a_{ii}< 0$ and $a_{ij}\geq 0$ for $i\neq j$, then all entries of $A^{-1}$ are nonpositive. Moreover, one can permute rows and columns so that $A^{-1}$ becomes block-diagonal with all entries in each block nonzero.\\ 

{\bf Definition (\cite{NikulinII}).} An at most hyperbolic system of vectors $\mathcal V$ is called:
\begin{itemize}
\item {\it elliptic}, if its Gram matrix $A$ is negative definite;
\item {\it connected parabolic}, if it is not elliptic, its Gram matrix $A$ is negative semidefinite and the associated weighted graph is connected;
\item {\it hyperbolic}, if its Gram matrix $A$ is not negative semidefinite;
\item {\it Lanner}, if it is hyperbolic, but any proper subsystem of it is not. 
\end{itemize}

{\bf Remarks (\cite{NikulinII}).} (1) A system is hyperbolic if and only if it contains a Lanner subsystem.\\

(2) The weighted graph of a Lanner system is connected. The subgraphs of the weighted graph of a hyperbolic system $\mathcal V$ corresponding to two Lanner subsystems $\mathcal V_1$ and $\mathcal V_2$ of $\mathcal V$ are connected by an edge.\\

(3) If $\mathcal V'$ is a proper subsystem of an elliptic or a connected parabolic system $\mathcal V$, then $\mathcal V'$ is elliptic. If $\mathcal V'$ is a proper subsystem of a Lanner system $\mathcal V$, then $\mathcal V'$ is either elliptic or connected parabolic (and in the latter case $\mathcal V'$ has one element less than $\mathcal V$).\\ 

(4) If $\mathcal V=\{ v_1, ..., v_n \}$ is an elliptic system, then its {\it canonical element}
$$
K(\mathcal V)={\alpha}_1 v_1 + ... + {\alpha}_n v_n, \;\;\; {\alpha}_i\in \mathbb Q
$$
exists and is uniquely defined by the conditions 
$$
K(\mathcal V)\cdot v_i=b_i-2, \; \;\; i=1,...,n.
$$

{\bf Definition.} (cf. \cite{NikulinII}) An at most hyperbolic system of vectors $\mathcal V$ is {\it log canonical}, if the canonical element of any elliptic subsystem $\mathcal V'=\{ v_1, ..., v_n \}$ of $\mathcal V$, $K(\mathcal V)=\sum_i {\alpha}_i v_i$, satisfies:
\begin{enumerate}
\item ${\alpha}_i\geq -1$ for any $i$, and
\item ${\alpha}_i> -1$, if the connected component of the vertex of $v_i$ in the graph associated to $\mathcal V'$ contains the vertex of an element $e$ with $e^2=-1$. 
\end{enumerate}

{\bf Remarks (\cite{NikulinII}).} (1) If $\mathcal V$ is a log canonical system, then its arbitrary subsystem $\mathcal V'\subset \mathcal V$ is also log canonical.\\

(2) If $\mathcal V$ is a log canonical system, then $v\in \mathcal V$ is {\it of the first kind} if $v^2=-1$, and $v\in \mathcal V$ is {\it of the second kind} if $v^2\leq -2$.\\

{\bf Definition (\cite{NikulinII}).} Let $\mathcal V =\{ v_1, ..., v_n, e \}$ be an at most hyperbolic system of vectors such that $e^2=-1$ and for any $i$ either $v_i\cdot e=0$, or $v_i^2\leq -2$ and $v_i\cdot e=1$.\\

Let $\widetilde{\mathcal V} =\{ \tilde v_1, ..., \tilde v_n \}$, where for any $i$
$$
\tilde v_i=v_i+(v_i\cdot e)e.
$$ 

We will call the passage from $\mathcal V$ to $\widetilde{\mathcal V}$ the {\it contraction} of $e$ and the passage from $\widetilde{\mathcal V}$ to $\mathcal V$ the {\it blow-up} of $\mathcal V'$, where $\mathcal V'\subset \widetilde{\mathcal V}$ is the subsystem consisting of all vectors $\tilde v_i$ such that $v_i\cdot e=1$.\\

{\bf Remark (\cite{NikulinII}).} $\widetilde{\mathcal V}$ is an at most hyperbolic system of vectors. If $\mathcal V$ is elliptic, connected parabolic, hyperbolic, Lanner or log canonical, then so is $\widetilde{\mathcal V}$.\\

{\bf Proposition A.1.1 (cf. \cite{NikulinII}, Proposition 1.1.1).} {\it The system of two vectors $\{  v_1, v_2 \}$ is elliptic and log canonical if and only if either
\begin{enumerate}
\item $v_1\cdot v_2=0$, $v_1^2< 0$, $v_2^2< 0$, or
\item $v_1\cdot v_2=1$, $v_1^2+v_2^2< -2$, $v_1^2< 0$, $v_2^2< 0$, or
\item $v_1\cdot v_2=2$, $v_1^2+v_2^2< -4$, $v_1^2< -1$, $v_2^2< -1$.
\end{enumerate}}

{\bf Corollary A.1.2 (\cite{NikulinII}, Corollary 1.1.2).} {\it If $\mathcal V$ is a log canonical system, then an element of the first kind $e\in \mathcal V$ is {\it contractible} (i.e. satisfies the assumptions of the definition of the contraction above) if and only if every system $\{ v, e \}$, $v\in \mathcal V$, $v\neq e$, is elliptic.}\\

{\bf Definition (\cite{NikulinII}).} A log canonical system is {\it minimal}, if it does not contain contractible elements of the first kind.\\

{\bf Theorem A.1.4 (cf. \cite{NikulinII}, Theorem~1.1.4).} {\it (a) An elliptic log canonical system is minimal if and only if it does not contain elements of the first kind. The connected associated weighted graphs of minimal elliptic log canonical systems are $\Gamma_1, ..., \Gamma_{6}$ (Picture~1), where all unspecified weights are assumed to be at least $2$ and the Gram matrix is negative definite.\\

(b) An elliptic system $\mathcal V$ is log canonical if and only if for every subsystem $\mathcal V' \subset \mathcal V$ with a connected associated graph $\Gamma '$ either
\begin{itemize}
\item $\Gamma '$ is one of $\Gamma_1, ..., \Gamma_{6}$ (Picture~1); or
\item $\mathcal V'$ is obtained by a sequence of blow-ups of subsystems with at most $2$ elements from a (possibly empty) minimal elliptic log terminal system in the sense of Nikulin \cite{NikulinII}, whose associated graph is either $\Gamma_1$ or $\Gamma_{2}$ (Picture~1).
\end{itemize}}

\begin{center}
\pspicture*[](-1,0)(14,8.5)
\rput(0,7.5){%
\psdots*[](0,0)(.5,0)(1,0)(2,0)
\qline(0,0)(1.4,0) \qline(1.6,0)(2,0)
\psline[linestyle=dotted]{}(1.5,.2)(1.5,-.2)
\rput(3,0){${\Gamma}_1$} \rput(1,-.5){$\underbrace{\qquad\qquad\quad}_{n\geq 0 \;\;\mbox{\small{vertices}}}$}
}
\rput(6,7){%
\psdots*[](0,0)(.5,0)(-2,0)(-1.5,0)(-.5,0)(1.5,0)(2,0)(0,-.5)(0,-1.5)(0,-2)
\qline(0,-1.1)(0,-2) \qline(0,0)(0,-.9) \qline(-2,0)(-1.1,0) \qline(-.9,0)(.9,0) \qline(1.1,0)(2,0)
\psline[linestyle=dotted]{}(-1,.2)(-1,-.2) \psline[linestyle=dotted]{}(1,.2)(1,-.2) \psline[linestyle=dotted]{}(.2,-1)(-.2,-1)
\rput(1.8,-2){${\Gamma}_2$} \rput(0,-3){$\frac{1}{p+1}+\frac{1}{q+1}+\frac{1}{r+1}\geq 1$} \rput(-1.25,-.5){$\underbrace{\qquad\qquad}_{p\geq 1 \;\;\mbox{\small{vertices}}}$} \rput(1.25,.5){$\overbrace{\qquad\qquad}^{q\geq 1 \;\;\mbox{\small{vertices}}}$} \rput(1.25,-1.25){$\left. \begin{tabular}{l} \\ \\ \\ \end{tabular}\right\} {{\scriptstyle r\geq 1 \; }\mbox{ \small{vertices}}}$}
}
\rput(10,7.3){%
\psdots*[](0,0)(0,.5)(0,-.5)(.5,0)(1,0)(2,0)(2,.5)(2,-.5)
\qline(0,.5)(0,-.5) \qline(2,.5)(2,-.5) \qline(0,0)(1.4,0) \qline(1.6,0)(2,0)
\psline[linestyle=dotted]{}(1.5,.2)(1.5,-.2)
\rput(3,0){${\Gamma}_3$} \rput(1,-.75){$\underbrace{\qquad\quad}_{n\geq 0 \;\;\mbox{\small{vertices}}}$} \rput(-.25,.5){$\scriptstyle{2}$} \rput(-.25,-.5){$\scriptstyle{2}$} \rput(2.25,.5){$\scriptstyle{2}$} \rput(2.25,-.5){$\scriptstyle{2}$}
}
\rput(.5,4.5){%
\psdots*[](0,0)(1,0)(-1,0)(0,1)(0,-1)
\qline(-1,0)(1,0) \qline(0,-1)(0,1)
\rput(2,-1){${\Gamma}_4$} \rput(1,.25){$\scriptstyle{2}$} \rput(-1,.25){$\scriptstyle{2}$} \rput(.25,1){$\scriptstyle{2}$} \rput(.25,-1){$\scriptstyle{2}$}
}
\rput(4.5,2){%
\psdots*[](0,0)(.5,1)(1.5,1)(2,0)(1.5,-1)(.5,-1)
\qline(0,0)(.5,1) \qline(.5,1)(1.5,1) \qline(1.5,1)(2,0) \qline(2,0)(1.5,-1) \qline(1.5,-1)(.5,-1) \qline(0,0)(0.125,-.25) \qline(.5,-1)(0.375,-.75)
\psline[linestyle=dotted]{}(.375,-.4)(.125,-.6)
\rput(4,-1){${\Gamma}_5$} \rput(4.25,0){$\left. \begin{tabular}{l} \\ \\ \\ \\ \\ \end{tabular}\right\} \mbox{\small{cycle with }} { \scriptstyle n\geq 3 \; }\mbox{ \small{vertices}} $}
}
\rput(10,4.5){%
\psdots*[](0,0)(1,0)
\qline(0,0.06)(1,0.06) \qline(0,-0.04)(1,-0.04)
\rput(2,0){${\Gamma}_6$}
}
\endpspicture

Picture $1$.\\
\end{center}

{\it Proof:} The proof is the same as that of \cite{NikulinII}, Theorem~1.1.4. One can use \cite{KM}, Theorem~4.7 (items 2, 3, 4 and 5) and \cite{KM}, Theorem~4.16 to list the graphs explicitly. {\it QED}\\

{\bf Corollary A.1.5 (cf. \cite{NikulinII}, Corollary 1.1.5).} {\it The associated graph $\Gamma$ of an elliptic log canonical system is a union of isolated cycles (with at least two vertices) and a tree without multiple edges. Every vertex of $\Gamma$ of weight $1$ has valency at most $2$ and every connected subgraph of $\Gamma$ consisting of vertices of weight $\geq 2$ is one of $\Gamma_1, ..., \Gamma_{6}$ (Picture~1).}\\

{\bf Theorem A.1.6 (cf. \cite{NikulinII}, Theorem~1.1.6).} {\it (a) A connected parabolic log canonical system is minimal if and only if either it does not contain elements of the first kind, or it has exactly two elements.\\

A minimal connected parabolic log canonical system has either the associated graph among $Q_1, ..., Q_{14}$ (Picture~3), or the Gram matrix of the form
$$
\left[
   \begin{array}{cc}
-b_1 & \sqrt{b_1b_2}\\
\sqrt{b_1b_2} & -b_2
   \end{array}
 \right], \;\;\; b_1,b_2 \geq 1.
$$

(b) Any connected parabolic log canonical system containing at least $3$ elements either has as its associated graph one of the graphs $Q_7, ..., Q_{14}$ (Picture~3), or is obtained by a sequence of blow-ups of subsystems with at most $2$ elements from a log canonical system, whose associated graph is either one of $Q_1, ..., Q_{6}$ (Picture~3) or one of $P_1$, ..., $P_8$ (Picture~2).}\\

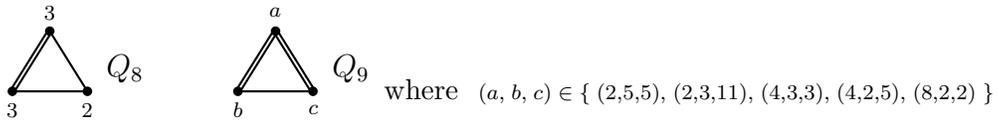
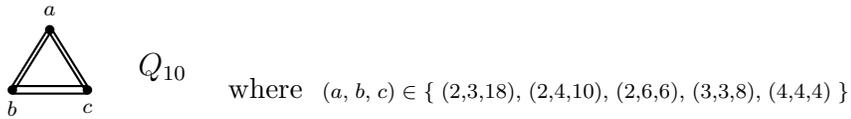
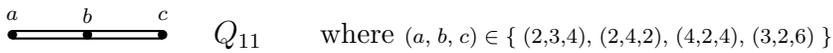
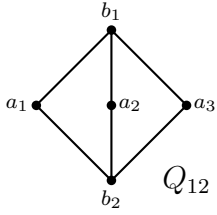
\begin{center}
\pspicture*[](-1,0)(14,4)
\rput(0,1){%
\psdots*[](1,0)(2,0)(3,0)(4,0)(0,2)(1,2)(2,2)(3,2)(4,2)(5,2)
\qline(0,2)(1,2) \qline(2,2.04)(3,2.04) \qline(2,1.96)(3,1.96) \qline(4,2.04)(5,2.04) \qline(4,1.96)(5,1.96) \qline(1,0.00)(2,0.00) \qline(1,0.05)(2,0.05) \qline(1,-0.05)(2,-0.05) \qline(3,0.00)(4,0.00) \qline(3,0.05)(4,0.05) \qline(3,-0.05)(4,-0.05)
\rput(.5,1.5){$P_1$} \rput(2.5,1.5){$P_2$} \rput(4.5,1.5){$P_3$} \rput(1.5,-.5){$P_7$} \rput(3.5,-.5){$P_8$}
\rput(-.2,2){${\scriptstyle 1}$} \rput(1.2,2){${\scriptstyle 1}$} \rput(1.8,2){${\scriptstyle 1}$} \rput(3.2,2){${\scriptstyle 4}$} \rput(3.8,2){${\scriptstyle 2}$} \rput(5.2,2){${\scriptstyle 2}$} \rput(.8,0){${\scriptstyle 3}$} \rput(2.2,0){${\scriptstyle 3}$} \rput(2.8,0){${\scriptstyle 1}$} \rput(4.2,0){${\scriptstyle 9}$}
}
\rput(7,2){%
\psdots*[](-.7,1)(.7,1)(0,0)(0,-1) \psdots*[](1.8,1)(2.5,0)(3.2,1)(2.5,-1) \psdots*[](4.3,1)(5,0)(5.7,1)(5,-1)
\qline(0,0)(-.7,1) \qline(0,0)(.7,1) \qline(0,0)(0,-1) \qline(2.5,0)(1.8,1) \qline(2.5,0)(3.2,1) \qline(2.5,0)(2.5,-1) \qline(5,0)(4.3,1) \qline(5,0)(5.7,1) \qline(5,0)(5,-1)
\rput(.7,-.7){$P_4$} \rput(3.2,-.7){$P_5$} \rput(5.7,-.7){$P_6$}
\rput(-.5,1){${\scriptstyle 3}$} \rput(.9,1){${\scriptstyle 6}$} \rput(.2,-.1){${\scriptstyle 1}$} \rput(0,-1.2){${\scriptstyle 2}$} \rput(2,1){${\scriptstyle 4}$} \rput(3.4,1){${\scriptstyle 4}$} \rput(2.7,-.1){${\scriptstyle 1}$} \rput(2.5,-1.2){${\scriptstyle 2}$} \rput(4.5,1){${\scriptstyle 3}$} \rput(5.9,1){${\scriptstyle 3}$} \rput(5.2,-.1){${\scriptstyle 1}$} \rput(5,-1.2){${\scriptstyle 3}$}
}
\endpspicture

Picture $2$.\\
\end{center}

\vspace{10mm}

\begin{center}
\pspicture*[](-1,0)(14,7)
\rput(.5,5){%
\psdots*[](0,0)(1,0)(-1,0)(0,1)(0,-1)
\qline(-1,0)(1,0) \qline(0,-1)(0,1)
\rput(1,-1){$Q_1$} \rput(1,.25){$\scriptstyle{2}$} \rput(-1,.25){$\scriptstyle{2}$} \rput(0,1.25){$\scriptstyle{2}$} \rput(0,-1.25){$\scriptstyle{2}$} \rput(.25,.25){$\scriptstyle{2}$}
}
\rput(4,5){%
\psdots*[](0,0)(0,1)(0,-1)(1,0)(2,0)(4,0)(4,1)(4,-1)
\qline(0,1)(0,-1) \qline(4,1)(4,-1) \qline(0,0)(2.9,0) \qline(3.1,0)(4,0)
\psline[linestyle=dotted]{}(3,.2)(3,-.2)
\rput(2,-1.3){$Q_2$} \rput(2,-.5){$\underbrace{\qquad\qquad\qquad\quad}_{n\geq 0 \;\;\mbox{\small{vertices}}}$} \rput(-.25,1){$\scriptstyle{2}$} \rput(-.25,0){$\scriptstyle{2}$} \rput(-.25,-1){$\scriptstyle{2}$} \rput(4.25,1){$\scriptstyle{2}$} \rput(4.25,0){$\scriptstyle{2}$} \rput(4.25,-1){$\scriptstyle{2}$} \rput(1,.25){$\scriptstyle{2}$} \rput(2,.25){$\scriptstyle{2}$}
}
\rput(10.5,5){%
\psdots*[](0,0)(1,0)(2,0)(3,0)(0,1)(1,1)(2,1)(0,-1)
\qline(0,-1)(0,1) \qline(0,1)(2,1) \qline(0,0)(3,0) 
\rput(2,-1){$Q_3$} \rput(-.25,1){$\scriptstyle{2}$} \rput(-.25,0){$\scriptstyle{2}$} \rput(-.25,-1){$\scriptstyle{2}$} \rput(1,1.25){$\scriptstyle{2}$} \rput(2,1.25){$\scriptstyle{2}$} \rput(1,.25){$\scriptstyle{2}$} \rput(2,.25){$\scriptstyle{2}$} \rput(3,.25){$\scriptstyle{2}$}
}
\rput(0,2.5){%
\psdots*[](0,0)(1,0)(2,0)(3,0)(4,0)(2,-1)(2,-2)
\qline(0,0)(4,0) \qline(2,0)(2,-2)
\rput(4,-1.5){$Q_4$} \rput(0,.25){$\scriptstyle{2}$} \rput(1,.25){$\scriptstyle{2}$} \rput(2,.25){$\scriptstyle{2}$} \rput(3,.25){$\scriptstyle{2}$} \rput(4,.25){$\scriptstyle{2}$} \rput(2.25,-1){$\scriptstyle{2}$} \rput(2.25,-2){$\scriptstyle{2}$}
}
\rput(6,2.5){%
\psdots*[](0,0)(1,0)(2,0)(3,0)(4,0)(5,0)(6,0)(7,0)(2,-1)
\qline(0,0)(7,0) \qline(2,0)(2,-1)
\rput(5,-1){$Q_5$} \rput(0,.25){$\scriptstyle{2}$} \rput(1,.25){$\scriptstyle{2}$} \rput(2,.25){$\scriptstyle{2}$} \rput(3,.25){$\scriptstyle{2}$} \rput(4,.25){$\scriptstyle{2}$} \rput(5,.25){$\scriptstyle{2}$} \rput(6,.25){$\scriptstyle{2}$} \rput(7,.25){$\scriptstyle{2}$} \rput(2.25,-1){$\scriptstyle{2}$}
}
\endpspicture
\end{center}

\begin{center}
\pspicture*[](-1,0)(14,5)
\rput(1,3){%
\psdots*[](0,0)(.5,1)(1.5,1)(2,0)(1.5,-1)(.5,-1)
\qline(0,0)(.5,1) \qline(.5,1)(1.5,1) \qline(1.5,1)(2,0) \qline(2,0)(1.5,-1) \qline(1.5,-1)(.5,-1) \qline(0,0)(0.125,-.25) \qline(.5,-1)(0.375,-.75)
\psline[linestyle=dotted]{}(.375,-.4)(.125,-.6)
\rput(3,-.5){$Q_6$} \rput(.5,1.25){$\scriptstyle{2}$} \rput(1.5,1.25){$\scriptstyle{2}$} \rput(-.25,0){$\scriptstyle{2}$} \rput(2.25,0){$\scriptstyle{2}$} \rput(.5,-1.25){$\scriptstyle{2}$} \rput(1.5,-1.25){$\scriptstyle{2}$} 
\rput(1,-2){$\underbrace{\qquad\qquad\qquad}_{\genfrac{}{}{0pt}{}{\mbox{\small cycle with } {\textstyle \; n\geq 3 \; } \mbox{\small{vertices}}}{ \mbox{\small{(all with weight}} {\textstyle \; 2)}}}$}
}
\rput(7.5,3){%
\psdots*[](-1,0)(0,0)(.5,1)(1.5,1)(2,0)(1.5,-1)(.5,-1)
\qline(0,0)(-1,0) \qline(0,0)(.5,1) \qline(.5,1)(1.5,1) \qline(1.5,1)(2,0) \qline(2,0)(1.5,-1) \qline(1.5,-1)(.5,-1) \qline(0,0)(0.5,-1)
\rput(-1.5,1){$Q_7:$} \rput(-1,.25){$\scriptstyle{2}$} \rput(-.15,.25){$\scriptstyle{2}$} \rput(2.25,0){$\scriptstyle{4}$} \rput(.5,1.25){$\scriptstyle{2}$} \rput(1.5,1.25){$\scriptstyle{2}$} \rput(.5,-1.25){$\scriptstyle{2}$} \rput(1.5,-1.25){$\scriptstyle{2}$} 
}
\rput(10,2){%
\psdots*[](0,0)(1.5,0)(3,0)(1.5,1.5)(3,1.5)
\qline(0,0)(3,0) \qline(1.5,0)(1.5,1.5) \qline(1.5,1.5)(3,1.5) \qline(3,1.5)(3,0) 
\rput(0,-.25){$\scriptstyle{2}$} \rput(1.5,-.25){$\scriptstyle{2}$} \rput(3,-.25){$\scriptstyle{2}$} \rput(1.5,1.75){$\scriptstyle{2}$} \rput(3,1.75){$\scriptstyle{3}$} 
}
\endpspicture

Picture $3$.\\
\end{center}

\begin{center}
\pspicture*[](-1,0)(14,8)
\rput(0,6){%
\psdots*[](0,0)(1,0)(0.5,0.8)
\qline(0,0)(1,0) \qline(1,0)(0.5,0.8) \qline(0,0.05)(0.5,0.85) \qline(0,-0.05)(0.5,0.75) 
\rput(1.5,.3){$Q_8$} \rput(0,-.25){$\scriptstyle{3}$} \rput(1,-.25){$\scriptstyle{2}$} \rput(0.5,1.05){$\scriptstyle{3}$}
}
\rput(3,6){%
\psdots*[](0,0)(1,0)(0.5,0.8)
\qline(0,0)(1,0) \qline(1,0.05)(0.5,0.85) \qline(1,-0.05)(0.5,0.75) \qline(0,0.05)(0.5,0.85) \qline(0,-0.05)(0.5,0.75)
\rput(1.5,.3){$Q_9$} \rput(0,-.25){$\scriptstyle{b}$} \rput(1,-.25){$\scriptstyle{c}$} \rput(0.5,1.05){$\scriptstyle{a}$} \rput(6,0){{\small where } ${\scriptstyle (a,\; b,\; c) \; \in \; \{ \; (2,5,5),\; (2,3,11),\; (4,3,3),\; (4,2,5),\; (8,2,2)\; \} }$}
}
\rput(0,3){%
\psdots*[](0,0)(1,0)(0.5,0.8)
\qline(0,0.05)(1,0.05) \qline(0,-0.05)(1,-0.05) \qline(1,0.05)(0.5,0.85) \qline(1,-0.05)(0.5,0.75) \qline(0,0.05)(0.5,0.85) \qline(0,-0.05)(0.5,0.75)
\rput(2,.3){$Q_{10}$} \rput(0,-.25){$\scriptstyle{b}$} \rput(1,-.25){$\scriptstyle{c}$} \rput(.5,1.05){$\scriptstyle{a}$} \rput(7,0){{\small where } ${\scriptstyle (a,\; b,\; c)\; \in \; \{ \; (2,3,18),\; (2,4,10),\; (2,6,6),\; (3,3,8),\; (4,4,4)\; \} }$}
}
\rput(0,.5){%
\psdots*[](0,0)(1,0)(2,0)
\qline(0,0.05)(2,0.05) \qline(0,-0.05)(2,-0.05)
\rput(3,0){$Q_{11}$} \rput(0,.25){$\scriptstyle{a}$} \rput(1,.25){$\scriptstyle{b}$} \rput(2,.25){$\scriptstyle{c}$} \rput(7.5,0){{\small where} ${\scriptstyle (a,\; b,\; c)\; \in \; \{ \; (2,3,4),\; (2,4,2),\; (4,2,4),\; (3,2,6)\; \} }$}
}
\endpspicture
\end{center}

\begin{center}
\pspicture*[](-1,0)(14,8)
\rput(1,6){%
\psdots*[](0,0)(1,0)(-1,0)(0,1)(0,-1)
\qline(0,-1)(0,1) \qline(0,-1)(-1,0) \qline(0,-1)(1,0) \qline(0,1)(-1,0) \qline(0,1)(1,0)
\rput(1,-1){$Q_{12}$} \rput(0,1.25){$\scriptstyle{b_1}$} \rput(0,-1.25){$\scriptstyle{b_2}$} \rput(-1.25,0){$\scriptstyle{a_1}$} \rput(1.25,0){$\scriptstyle{a_3}$} \rput(.25,0){$\scriptstyle{a_2}$}
\rput(7,0){ $\genfrac{}{}{0pt}{}{\mbox{\small where }\;(b_1,\; b_2,\; a_1,\; a_2,\; a_3)\; \in \; \{\; (2,2,2,3,6),\; (2,2,2,4,4),\; (2,2,3,3,3),\; (2,3,2,2,5),}{ (2,4,2,2,3),\; (2,6,2,2,2),\; (3,3,2,2,2)\; \} }$}
}
\rput(-.5,1.5){%
\psdots*[](0,0)(2,0)(1,0.7)(1,1.5)
\qline(0,0)(2,0) \qline(0,0)(1,0.7) \qline(0,0)(1,1.5) \qline(2,0)(1,0.7) \qline(2,0)(1,1.5) \qline(1,0.7)(1,1.5)
\rput(1,-1){$Q_{13}$} \rput(1,1.75){$\scriptstyle{a_1}$} \rput(1,0.45){$\scriptstyle{a_4}$} \rput(0,-.25){$\scriptstyle{a_2}$} \rput(2,-.25){$\scriptstyle{a_3}$}
\rput(4,1){$\genfrac{}{}{0pt}{}{\mbox{\small where }\;(a_1,\; a_2,\; a_3,\; a_4)\; \in \; \{\; (2,2,5,5), }{ (2,2,3,11),\; (2,3,3,5),\; (3,3,3,3) \; \} }$}
}
\rput(7.75,1.5){%
\psdots*[](0,0)(-1,.5)(1,.5)(0,1)
\qline(0,0)(-1,.5) \qline(0,0)(1,.5) \qline(0,0)(0,1) \qline(0,1)(-1,.5) \qline(0,1)(1,.5) 
\rput(0,-1){$Q_{14}$} \rput(0,-.25){$\scriptstyle{a_2}$} \rput(0,1.25){$\scriptstyle{a_1}$} \rput(-1.25,.5){$\scriptstyle{b_2}$} \rput(1.25,.5){$\scriptstyle{b_1}$}
\rput(3.75,0){$\genfrac{}{}{0pt}{}{\mbox{\small where }\;(a_1,\; a_2,\; b_1,\; b_2)\; \in \; \{ \; (2,2,4,4), }{ (2,2,3,6),\; (2,5,2,2),\; (3,3,2,2)\; \} }$}
}

\endpspicture

Picture $3$ (continuation).\\
\end{center}

{\it Proof:} The argument is essentially the same as in \cite{NikulinII}, Theorem~1.1.6.\\

Let $\Gamma$ be the associated graph of a minimal connected parabolic log canonical system with at least $3$ elements. If $\Gamma$ contains a double edge, then it has exactly $3$ elements (since any minimal elliptic log canonical system, whose associated graph is connected and has a double edge, has exactly $2$ elements). Then one checks that $\Gamma$ is one of $Q_8$, $Q_9$, $Q_{10}$, $Q_{11}$.\\

Note that neither of $Q_8$, $Q_9$, $Q_{10}$, $Q_{11}$ can appear as a result of a contraction of an element of the first kind of a connected parabolic log canonical system.\\

Assume that $\Gamma$ does not contain a double edge. If $\Gamma$ contains a cycle, then it is one of $Q_6$, $Q_7$, $Q_{12}$, $Q_{13}$, $Q_{14}$.\\

Note that one can not obtain a connected parabolic log canonical system by blowing up a subsystem of $Q_{7}$, $Q_{12}$, $Q_{13}$, $Q_{14}$.\\

Hence we may assume that $\Gamma$ is a tree without double edges. So, it is obtained from ${\Gamma}_1$, ${\Gamma}_2$, ${\Gamma}_3$ or ${\Gamma}_4$ by adding a single vertex and an edge. We get $Q_1$, $Q_2$, $Q_3$, $Q_4$, $Q_5$.\\

Note that there is no minimal connected parabolic log canonical system with the associated graph of the following form:

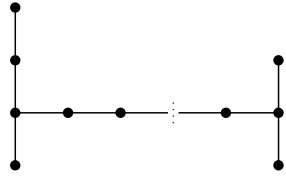
\begin{center}
\psscalebox{.7}{%
\pspicture*[](0,0)(7,5)
\psdots*[dotsize=4pt 2](1,1)(1,2)(1,3)(1,4)(2,2)(3,2)(5,2)(6,2)(6,1)(6,3)
\qline(1,1)(1,4)
\qline(1,2)(3.9,2)
\qline(4.1,2)(6,2)
\qline(6,3)(6,1)
\psline[linestyle=dotted]{}(4,2.2)(4,1.8)
\endpspicture
}
\end{center}

{\it QED}\\

{\bf Corollary A.1.7 (cf. \cite{NikulinII}, Corollary 1.1.7).} {\it If $\Gamma$ is the associated graph of a connected parabolic log canonical system with at least $10$ elements, then $\Gamma$ is either a tree or a cycle. It has no multiple edges and every vertex of $\Gamma$ with weight $1$ has valency at most $2$.}\\

{\bf Theorem A.1.8 (cf. \cite{NikulinII}, Theorem~1.1.8).} {\it (a) A Lanner log canonical system is minimal if and only if either it does not contain elements of the first kind, or it has exactly $2$ elements, or it has exactly $3$ elements and its Gram matrix is different (upto an equivalence) from the following ones:
$$
\left[
   \begin{array}{ccc}
-1 & 1 & 0\\
1 & -b_1 & r\\
0 & r & -b_2
   \end{array}
 \right], \;\;\; r\geq 1, \;  b_1\geq 2, \; b_2 \geq 1;
$$

$$
\left[
   \begin{array}{ccc}
-1 & 1 & 1\\
1 & -b_1 & r\\
1 & r & -b_2
   \end{array}
 \right], \;\;\; r\geq 1, \;  b_1\geq 2, \; b_2 \geq 2.
$$

The associated graph of a minimal Lanner log canonical system with at least $4$ elements is one of the graphs $H_{8,1}$, ..., $H_{11}$ (Picture~4), where all unspecified weights are equal to $2$ and $b\geq 2$. If a minimal Lanner log canonical system has at most $3$ elements, then its Gram matrix is either one of the following:
$$
G_1(b_1,b_2;r)=\left[
   \begin{array}{cc}
-b_1 & r\\
r & -b_2
   \end{array}
 \right], \;\;\; r^2>b_1b_2;
$$

$$
G_2(b_1,b_2,b_3)=\left[
   \begin{array}{ccc}
-b_1 & \sqrt{b_1b_2} & 0\\
\sqrt{b_1b_2} & -b_2 & 1\\
0 & 1 & -b_3
   \end{array}
 \right], \;\;\; b_3 \geq 2;
$$

$$
G_3(b_1,b_2,b_3)=\left[
   \begin{array}{ccc}
-b_1 & \sqrt{b_1b_2} & 1\\
\sqrt{b_1b_2} & -b_2 & 1\\
1 & 1 & -b_3
   \end{array}
 \right], \;\;\; b_3 \geq 2;
$$

$$
G_4(b_1,b_2,b_3)=\left[
   \begin{array}{ccc}
-b_1 & \sqrt{b_1b_2} & 0\\
\sqrt{b_1b_2} & -b_2 & \sqrt{b_2b_3}\\
0 & \sqrt{b_2b_3} & -b_3
   \end{array}
 \right];
$$

$$
G_5(b_1,b_2,b_3)=\left[
   \begin{array}{ccc}
-b_1 & \sqrt{b_1b_2} & 1\\
\sqrt{b_1b_2} & -b_2 & \sqrt{b_2b_3}\\
1 & \sqrt{b_2b_3} & -b_3
   \end{array}
 \right], \;\;\; b_1+b_3 \geq 3;
$$

$$
G_6(b_1,b_2,b_3)=\left[
   \begin{array}{ccc}
-b_1 & \sqrt{b_1b_2} & \sqrt{b_1b_3}\\
\sqrt{b_1b_2} & -b_2 & \sqrt{b_2b_3}\\
\sqrt{b_1b_3} & \sqrt{b_2b_3} & -b_3
   \end{array}
 \right];
$$

$$
G_7(b_1,b_2,b_3)=\left[
   \begin{array}{ccc}
-b_1 & \sqrt{b_1b_2} & 0\\
\sqrt{b_1b_2} & -b_2 & 2\\
0 & 2 & -b_3
   \end{array}
 \right], \;\;\; b_2,b_3 \geq 2,\; b_2+b_3 \geq 5;
$$

$$
G_8(b_1,b_2,b_3)=\left[
   \begin{array}{ccc}
-b_1 & \sqrt{b_1b_2} & 1\\
\sqrt{b_1b_2} & -b_2 & 2\\
1 & 2 & -b_3
   \end{array}
 \right], \;\;\; b_2,b_3 \geq 2,\; b_2+b_3 \geq 5;
$$

$$
G_9(b_1,b_2,b_3)=\left[
   \begin{array}{ccc}
-b_1 & \sqrt{b_1b_2} & 2\\
\sqrt{b_1b_2} & -b_2 & 2\\
2 & 2 & -b_3
   \end{array}
 \right], \;\;\; b_1,b_2,b_3 \geq 2, \; b_1+b_3, b_2+b_3 \geq 5;
$$

$$
G_{10}(b_1,b_2,b_3)=\left[
   \begin{array}{ccc}
-b_1 & \sqrt{b_1b_2} & 2\\
\sqrt{b_1b_2} & -b_2 & \sqrt{b_2b_3}\\
2 & \sqrt{b_2b_3} & -b_3
   \end{array}
 \right], \;\;\; b_1,b_3 \geq 2, \; b_1+b_3 \geq 5;
$$

where $b_1, b_2, b_3\geq 1$ everywhere, or one of the following:

\begin{multline*}
G_{11}(b_1,b_2,b_3)=\left[
   \begin{array}{ccc}
-b_1 & 2 & 0\\
2 & -b_2 & 2\\
0 & 2 & -b_3
   \end{array}
 \right], \;\;\; (b_1,b_2,b_3) \in \{ (2,3,2), (2,3,3), (3,2,3), (3,2,4),\\
 (3,2,5) \};
\end{multline*}

$$
G_{12}(b_1,b_2,b_3)=\left[
   \begin{array}{ccc}
-b_1 & 1 & 1\\
1 & -b_2 & 2\\
1 & 2 & -b_3
   \end{array}
 \right], \;\;\; (b_1,b_2,b_3) \in \{ (2,2,3), (2,2,4), (3,2,3), (4,2,3) \};
$$

\begin{multline*}
G_{13}(b_1,b_2,b_3)=\left[
   \begin{array}{ccc}
-b_1 & 2 & 2\\
2 & -b_2 & 1\\
2 & 1 & -b_3
   \end{array}
 \right], \;\;\; (b_1,b_2,b_3) \in \{ (3,3,3), (3,4,3) \}\;\; \mbox{or}\\
 b_1=3,\; b_2=2,\; 2\leq b_3\leq 9\;\; \mbox{or}\;\; b_1=2,\; b_2=3,\; 3\leq b_3\leq 10\;\; \mbox{or}\\
  b_1=2,\; b_2=4,\; 4\leq b_3\leq 6;
\end{multline*}

\begin{multline*}
G_{14}(b_1,b_2,b_3)=\left[
   \begin{array}{ccc}
-b_1 & 2 & 2\\
2 & -b_2 & 2\\
2 & 2 & -b_3
   \end{array}
 \right], \;\;\; b_1=2,\; b_2=3,\; 3\leq b_3\leq 17\;\; \mbox{or}\\
 b_1=2,\; b_2=4,\; 4\leq b_3\leq 9\;\; \mbox{or}\;\; b_1=2,\; b_2=5,\; 5\leq b_3\leq 7\;\; \mbox{or}\\
 b_1=3,\; b_2=3,\; 3\leq b_3\leq 7\;\; \mbox{or}\;\; b_1=3,\; b_2=4,\; 4\leq b_3\leq 5.
\end{multline*}

(b) Any Lanner log canonical system containing at least $4$ elements either has the associated graph among $H_{8,1}$, ..., $H_{0,14,1}$, $H_{0,4}$, $H_{1}$, ..., $H_4$ (Picture~4), or is obtained by a sequence of blow-ups from a Lanner log canonical system with the associated graph among $H_{0,1}$, $H_{0,2}$, $H_{0,3}$, $H_0$, $H_5$, ..., $H_{11}$, or is obtained by a sequence of blow-ups from a Lanner log canonical system with the Gram matrix \\

$G_1(1,1;2)$, $G_1(1,2;2)$, $G_1(1,3;2)$, $G_1(1,1;3)$, $G_1(2,2;3)$, $G_1(2,2;4)$, $G_2(1,1,b)$, or\\

$G_2(2,2,b)$, $G_2(1,4,b)$, $G_2(4,1,b)$, $G_3(1,1,b)$, $G_3(2,2,b)$, $G_3(1,4,b)$, $G_3(1,9,b)$, or\\

$G_3(3,3,b)$, $G_4(1,1,1)$, $G_4(1,1,4)$, $G_4(4,1,4)$, $G_5(4,1,1)$, $G_5(4,1,4)$, $G_5(2,2,2)$, or\\

$G_5(9,1,1)$, $G_6(1,1,1)$, $G_6(1,1,4)$, $G_6(2,2,2)$, $G_6(3,3,3)$, $G_7(1,4,2)$, $G_7(2,2,4)$, or\\

$G_7(2,2,6)$, $G_7(4,1,4)$, $G_8(1,4,b)$, $G_8(2,2,b)$, $G_9(3,3,4)$, $G_9(2,2,b)$.}\\

\begin{center}
\pspicture*[](-1,0)(14,7.5)
\rput(0,5.5){%
\psdots*[](0,0)(1,0)(2,0)(0,1)
\qline(0,0)(2,0) \qline(0,1)(1,0) \qline(-.05,0)(-.05,1) \qline(0.05,0)(0.05,1)
\rput(1,-1){$H_{8,1}\;\; (b\geq 2)$} \rput(-.25,0){$\scriptstyle{3}$} \rput(-.25,1){$\scriptstyle{3}$} \rput(1,-.25){$\scriptstyle{2}$} \rput(2,-.25){$\scriptstyle{b}$}
}
\rput(4,5.5){%
\psdots*[](0,0)(1.5,0)(1.5,1.5)(0,1.5)
\qline(0,0)(1.5,0) \qline(0,0)(0,1.5) \qline(1.5,1.5)(1.5,0) \qline(1.5,1.5)(0,1.5) \qline(0,0.07)(1.5,1.57) \qline(0,-0.07)(1.5,1.43)
\rput(.7,-1){$H_{8,2}$} \rput(-.25,0){$\scriptstyle{3}$} \rput(-.25,1.5){$\scriptstyle{2}$} \rput(1.75,0){$\scriptstyle{2}$} \rput(1.75,1.5){$\scriptstyle{3}$}
}
\rput(7,5.5){%
\psdots*[](0,0)(2,0)(1,0.7)(1,1.5)
\qline(0,0)(2,0) \qline(0,0)(1,0.7) \qline(0,0.07)(1,1.57) \qline(0,-0.07)(1,1.43) \qline(2,0)(1,0.7) \qline(2,0)(1,1.5) \qline(1,0.7)(1,1.5)
\rput(1,-1){$H_{8,3}$} \rput(1,1.75){$\scriptstyle{3}$} \rput(1,0.45){$\scriptstyle{2}$} \rput(0,-.25){$\scriptstyle{3}$} \rput(2,-.25){$\scriptstyle{2}$}
}
\rput(10,5.5){%
\psdots*[](0,0)(1.5,0)(1.5,1.5)(0,1.5)
\qline(0,0.05)(1.5,0.05) \qline(0,-.05)(1.5,-.05) \qline(0.05,0)(0.05,1.5) \qline(-.05,0)(-.05,1.5) \qline(1.5,1.5)(1.5,0) \qline(1.5,1.5)(0,1.5) \qline(0,0.07)(1.5,1.57) \qline(0,-0.07)(1.5,1.43)
\rput(.7,-1){$H_{9,1}$} \rput(-.25,0){$\scriptstyle{4}$} \rput(-.25,1.5){$\scriptstyle{2}$} \rput(1.75,0){$\scriptstyle{2}$} \rput(1.75,1.5){$\scriptstyle{5}$}
}
\rput(0,1.5){%
\psdots*[](0,0)(2,0)(1,0.7)(1,1.5)
\qline(0,0)(2,0) \qline(0,0)(1,0.7) \qline(0,0.07)(1,1.57) \qline(0,-0.07)(1,1.43) \qline(2,0)(1,0.7) \qline(2,0.07)(1,1.57) \qline(2,-.07)(1,1.43) \qline(1.05,0.7)(1.05,1.5) \qline(.95,0.7)(.95,1.5)
\rput(1,-1){$H_{9,2}$} \rput(1,1.75){$\scriptstyle{a}$} \rput(1,0.45){$\scriptstyle{b}$} \rput(0,-.25){$\scriptstyle{b}$} \rput(2,-.25){$\scriptstyle{b}$} \rput(4.5,.7){{\small where} ${\scriptstyle (a,\;b)\;\in \; \{ \; (2,5), \; (4,3),\; (8,2)\; \} }$}
}
\rput(9,1.5){%
\psdots*[](0,0)(1,0)(2,0)(0,1)
\qline(0,0.05)(2,0.05) \qline(0,-.05)(2,-.05) \qline(0,1.07)(1,0.07) \qline(0,.93)(1,-.07) \qline(-.05,0)(-.05,1) \qline(0.05,0)(0.05,1)
\rput(2.5,.5){$H_{10,1}$} \rput(-.25,0){$\scriptstyle{6}$} \rput(-.25,1){$\scriptstyle{6}$} \rput(1,-.25){$\scriptstyle{2}$} \rput(2,-.25){$\scriptstyle{3}$}
}
\endpspicture

Picture $4$.\\
\end{center}

\begin{center}
\pspicture*[](-1,0)(14,7)
\rput(0,5){%
\psdots*[](0,0)(2,0)(1,0.7)(1,1.5)
\qline(0,0.05)(2,0.05) \qline(0,-0.05)(2,-0.05) \qline(0,.05)(1,0.77) \qline(0,-.05)(1,0.63) \qline(0,0.07)(1,1.57) \qline(0,-0.07)(1,1.43) \qline(2,0.07)(1,0.77) \qline(2,-.07)(1,0.63) \qline(2,0.07)(1,1.57) \qline(2,-.07)(1,1.43) \qline(1.05,0.7)(1.05,1.5) \qline(.95,0.7)(.95,1.5)
\rput(1,-1){$H_{10,2}$} \rput(1,1.75){$\scriptstyle{4}$} \rput(1,0.45){$\scriptstyle{4}$} \rput(0,-.25){$\scriptstyle{4}$} \rput(2,-.25){$\scriptstyle{4}$}
}
\rput(3,5.5){%
\psdots*[](0,0)(1,0)(2,0)(3,0) \psdots*[](0,1)(1,1)(2,1)(3,1)
\qline(0,0.05)(3,0.05) \qline(0,-0.05)(3,-0.05) \qline(0,1.05)(3,1.05) \qline(0,0.95)(3,0.95)
\rput(1.5,-1){$H_{11,1}$} \rput(0,.25){$\scriptstyle{2}$} \rput(1,.25){$\scriptstyle{4}$} \rput(2,.25){$\scriptstyle{2}$} \rput(3,.25){$\scriptstyle{4}$} \rput(0,1.25){$\scriptstyle{4}$} \rput(1,1.25){$\scriptstyle{3}$} \rput(2,1.25){$\scriptstyle{2}$} \rput(3,1.25){$\scriptstyle{6}$}
}
\rput(7,5.5){%
\psdots*[](0,0)(1,0)(2,0)(1,1)
\qline(0,0.05)(2,0.05) \qline(0,-.05)(2,-.05) \qline(1.05,1)(1.05,0) \qline(.95,1)(.95,0)
\rput(3.5,.5){$H_{11,2} \; \; {\scriptstyle ( \;a \; \in \; \{\; 2,\; 4\; \}\;  ) }$} \rput(0,-.25){$\scriptstyle{a}$} \rput(1,-.25){$\scriptstyle{6-a}$} \rput(2,-.25){$\scriptstyle{a}$} \rput(1,1.25){$\scriptstyle{a}$}
}
\rput(0,1.5){%
\psdots*[](0,0)(1.5,0)(1.5,1.5)(0,1.5)
\qline(0,0.05)(1.5,0.05) \qline(0,-.05)(1.5,-.05) \qline(0.05,0)(0.05,1.5) \qline(-.05,0)(-.05,1.5) \qline(1.55,1.5)(1.55,0) \qline(1.45,1.5)(1.45,0) \qline(1.5,1.55)(0,1.55) \qline(1.5,1.45)(0,1.45)
\rput(.7,-1){$H_{11,3}$} \rput(-.25,0){$\scriptstyle{4}$} \rput(-.25,1.5){$\scriptstyle{2}$} \rput(1.75,0){$\scriptstyle{2}$} \rput(1.75,1.5){$\scriptstyle{4}$}
}
\rput(3,1.5){%
\psdots*[](1,0)(0,1)(1,1)(2,1)
\qline(0,1.05)(2,1.05) \qline(0,.95)(2,.95) \qline(1.05,1)(1.05,0) \qline(.95,1)(.95,0) \qline(0,1)(1,0) \qline(2,1)(1,0)
\rput(1,-1){$H_{11,4}$} \rput(0,1.25){$\scriptstyle{2}$} \rput(1,-.25){$\scriptstyle{5}$} \rput(1,1.25){$\scriptstyle{4}$} \rput(2,1.25){$\scriptstyle{2}$}
}
\rput(7,1.5){%
\psdots*[](1,0)(0,1)(1,1)(2,1)
\qline(0,1.05)(2,1.05) \qline(0,.95)(2,.95) \qline(1,1)(1,0) \qline(0,1.07)(1,0.07) \qline(0,.93)(1,-.07) \qline(2,1.07)(1,0.07) \qline(2,.93)(1,-.07)
\rput(1,-1){$H_{11,5}$} \rput(0,1.25){$\scriptstyle{4}$} \rput(1,-.25){$\scriptstyle{5}$} \rput(1,1.25){$\scriptstyle{2}$} \rput(2,1.25){$\scriptstyle{4}$}
}
\rput(10,1.5){%
\psdots*[](1,0)(0,1)(1,1)(2,1)
\qline(0,1.05)(2,1.05) \qline(0,.95)(2,.95) \qline(1.05,1)(1.05,0) \qline(.95,1)(.95,0) \qline(0,1.07)(1,0.07) \qline(0,.93)(1,-.07) \qline(2,1.07)(1,0.07) \qline(2,.93)(1,-.07)
\rput(1,-1){$H_{11,6} \;\; {\scriptstyle (\; a \;\in \; \{ \; 2,\; 4 \; \} \; ) }$} \rput(0,1.25){$\scriptstyle{a}$} \rput(1,-.25){$\scriptstyle{10}$} \rput(1,1.25){$\scriptstyle{6-a}$}
\rput(2,1.25){$\scriptstyle{a}$}
}
\endpspicture
\end{center}

\begin{center}
\pspicture*[](-1,0)(14,8)
\rput(-.5,6){%
\psdots*[](1,0)(0,0)(2,0)(2.5,1)(3.5,1)(4,0)(2.5,-1)(3.5,-1)
\qline(0,0)(2,0) \qline(2,0)(2.5,1) \qline(2.5,1)(3.5,1) \qline(3.5,1)(4,0) \qline(4,0)(3.5,-1) \qline(3.5,-1)(2.5,-1) \qline(2,0)(2.5,-1)
\rput(1,-1){$H_{0,7,1}$} \rput(4.25,0){$\scriptstyle{4}$}
}
\rput(6,6){%
\psdots*[](0,0)(1.5,0)(-1,0)(1.5,1.5)(0,1.5)(0,-1)
\qline(-1,0)(1.5,0) \qline(0,-1)(0,1.5) \qline(0,1.5)(1.5,1.5) \qline(1.5,1.5)(1.5,0) 
\rput(1.75,1.5){$\scriptstyle{3}$} \rput(1,-1){$H_{0,7,2}$}
}
\rput(10,6){%
\psdots*[](0,0)(-1,0)(1,0)(0,1)(0,-1)(-1,1.5)
\qline(0,-1)(-1,0) \qline(0,-1)(1,0) \qline(0,1)(-1,0) \qline(0,1)(1,0) \qline(-1,0)(-1,1.5) \qline(0,-1)(0,1) \qline(-0.45,0.65)(-1,1.5) \qline(-0.35,0.55)(0,0)
\pscurve[fillstyle=none]{}(-1,1.5)(0,1.4)(1.2,1)(1,0)
\rput(0,-1.25){$\scriptstyle{3}$} \rput(-.15,1.15){$\scriptstyle{3}$} \rput(-1.25,1.5){$\scriptstyle{3}$} \rput(2,0){$H_{0,12,1}$}
}
\rput(.5,2.5){%
\psdots*[](0,0)(-.5,0)(.5,0)(0,1)(0,-1)(-1,0)
\qline(0,-1)(-1,0) \qline(0,-1)(-.5,0) \qline(0,-1)(.5,0) \qline(0,-1)(0,1) \qline(0,1)(-1,0) \qline(0,1)(-.5,0) \qline(0,1)(.5,0)
\rput(0,-1.25){$\scriptstyle{3}$} \rput(0,1.25){$\scriptstyle{3}$} \rput(0,-2){$H_{0,12,2}$}
}
\rput(3.5,1.5){%
\psdots*[](0,0)(1,.5)(-1,.5)(0,1)(0,2)
\qline(1,.5)(-1,.5) \qline(1,.5)(0,0) \qline(1,.5)(0,1) \qline(1,.5)(0,2) \qline(-1,.5)(0,0) \qline(-1,.5)(0,1) \qline(-1,.5)(0,2) \qline(0,2)(0,0.6) \qline(0,.4)(0,0)
\pscurve[]{}(0,2)(1.5,1.5)(1.4,0)(1,-.3)(0,0)
\rput(0,-.25){$\scriptstyle{3}$} \rput(-1.25,.5){$\scriptstyle{3}$} \rput(1.25,.5){$\scriptstyle{3}$} \rput(0,2.25){$\scriptstyle{3}$} \rput(-.15,1.15){$\scriptstyle{3}$}
\rput(0,-1){$H_{0,13,1}$}
}
\rput(8,2){%
\psdots*[](-.5,0)(.5,0)(0,1)(0,-1)(-1,0)
\qline(0,-1)(-1,0) \qline(0,-1)(-.5,0) \qline(0,-1)(.5,0) \qline(0,-1)(0,1) \qline(0,1)(-1,0) \qline(0,1)(-.5,0) \qline(0,1)(.5,0)
\rput(0,-1.25){$\scriptstyle{b}$} \rput(0,1.25){$\scriptstyle{a}$} \rput(3,0){{\small where } ${\scriptstyle (a,\; b)\; \in \; \{ \; (2,5),\; (3,3) \; \}  }$} \rput(-1,-1){$H_{0,14,1}$}
}
\endpspicture

Picture $4$ (continuation).\\
\end{center}

\begin{center}
\pspicture*[](-1,0)(14,8)
\rput(.5,5.5){%
\psdots*[](0,0)(1,0)(-1,0)(0,1)(0,-1)
\qline(0,-1)(0,1) \qline(0,-1)(-1,0) \qline(0,-1)(1,0) \qline(0,1)(-1,0) \qline(0,1)(1,0)
\rput(1.5,-1){$H_{0,1}$} \rput(0,1.25){$\scriptstyle{a}$} \rput(0,-1.25){$\scriptstyle{2}$} \rput(-1.25,0){$\scriptstyle{b_1}$} \rput(1.25,0){$\scriptstyle{2}$} \rput(.25,0){$\scriptstyle{b_2}$}
\rput(5,0.5){ $\genfrac{}{}{0pt}{}{\mbox{\small where }\;(a,\; b_1,\; b_2)\; \in \; \{ \; (2,2,x),\; (2,3,3),\; (2,3,4),\; (2,3,5),\; (3,2,2),}{(3,2,3),\; (3,2,4),\; (4,2,2),\; (5,2,2)\; \} , \;\; x\geq 2 }$}
}
\rput(4.5,3.5){%
\psdots*[](0,0)(2,0)(1,0.7)(1,1.5)
\qline(0,0)(2,0) \qline(0,0)(1,0.7) \qline(0,0)(1,1.5) \qline(2,0)(1,0.7) \qline(2,0)(1,1.5) \qline(1,0.7)(1,1.5)
\rput(1,-.5){$H_{0,2}$} \rput(1,1.75){$\scriptstyle{a}$} \rput(1,0.45){$\scriptstyle{2}$} \rput(0,-.25){$\scriptstyle{b}$} \rput(2,-.25){$\scriptstyle{c}$}
\rput(5.5,.5){$\genfrac{}{}{0pt}{}{\mbox{\small where }\;(a,\; b,\; c)\;\in \; \{\; (2,2,x),\; (2,3,y),\; (2,4,4),\; (2,4,5),}{(2,4,6),\; (3,3,3),\; (3,3,4)\; \} ,\;\; x\geq 2 , \;\; 3\leq y\leq 10 }$}
}
\rput(.5,1.5){%
\psdots*[](0,0)(-1,.5)(1,.5)(0,1)
\qline(0,0)(-1,.5) \qline(0,0)(1,.5) \qline(0,0)(0,1) \qline(0,1)(-1,.5) \qline(0,1)(1,.5) 
\rput(0,-1){$H_{0,3}$} \rput(0,-.25){$\scriptstyle{a}$} \rput(0,1.25){$\scriptstyle{2}$} \rput(-1.25,.5){$\scriptstyle{b}$} \rput(1.25,.5){$\scriptstyle{c}$}
\rput(5,0){$\genfrac{}{}{0pt}{}{\mbox{\small where }\;(a,\; b,\; c)\;\in \; \{ \; (2,2,x),\; (2,3,y),\; (3,2,2),\;  (3,2,3),}{(3,2,4),\; (4,2,2)\; \} ,\;\; x\geq 2 ,\;\; 3\leq y\leq 5 } $}
}
\rput(11,1.5){%
\psdots*[](0,0)(-1,0)(1,0)(0,1)(0,-1)
\qline(-1,0)(1,0) \qline(0,-1)(0,1) \qline(0,-1)(-1,0) \qline(0,-1)(1,0) \qline(0,1)(-1,0) \qline(0,1)(1,0) 
\rput(1.5,-1){$H_{0,4}$} \rput(-0.15,.15){$\scriptstyle{5}$}
}
\endpspicture
\end{center}

\begin{center}
\pspicture*[](-1,0)(14,8)
\rput(.5,6){%
\psdots*[](0,0)(.5,1)(1.5,1)(2,0)(1.5,-1)(.5,-1)(-1,0)
\qline(0,0)(.5,1) \qline(.5,1)(1.5,1) \qline(1.5,1)(2,0) \qline(2,0)(1.5,-1) \qline(1.5,-1)(.5,-1) \qline(0,0)(0.125,-.25) \qline(.5,-1)(0.375,-.75) \qline(-1,0)(0,0)
\psline[linestyle=dotted]{}(.375,-.4)(.125,-.6)
\rput(-1,1.5){$H_0:$}
\rput(1,-2){$\underbrace{\qquad\qquad\qquad}_{\genfrac{}{}{0pt}{}{\mbox{\small{cycle with }} {\textstyle \; n \; } \mbox{\small{vertices}},}{\displaystyle 3\leq n \leq 8}}$}
}
\rput(5,6){%
\psdots*[](-.5,0)(0,0)(.5,1)(1.5,.5)(1.5,-.5)(.5,-1)
\qline(0,0)(-.5,0) \qline(0,0)(.5,1) \qline(.5,1)(1.5,.5) \qline(1.5,.5)(1.5,-.5) \qline(1.5,-.5)(.5,-1) \qline(0,0)(0.5,-1)
\rput(-.5,.25){$\scriptstyle{b}$} \rput(.5,-2){${(3\leq b\leq 6)}$}
}
\rput(8,6){%
\psdots*[](-.5,0)(0,0)(.75,1)(1.5,0)(.75,-1)
\qline(0,0)(-.5,0) \qline(0,0)(.75,1) \qline(.75,1)(1.5,0) \qline(1.5,0)(.75,-1) \qline(0,0)(0.75,-1)
\rput(-.5,.25){$\scriptstyle{b}$} \rput(.5,-2){${(b\geq 3)}$}
}
\rput(11,6){%
\psdots*[](-.5,0)(0,0)(1,1)(1,-1)
\qline(0,0)(-.5,0) \qline(0,0)(1,1) \qline(1,1)(1,-1) \qline(0,0)(1,-1)
\rput(-.5,.25){$\scriptstyle{b}$} \rput(0,-2){${(b\geq 3)}$}
}
\rput(0,1.5){%
\psdots*[](-.5,0)(0,0)(.5,1)(1.5,.5)(1.5,-.5)(.5,-1)
\qline(0,0)(-.5,0) \qline(0,0)(.5,1) \qline(.5,1)(1.5,.5) \qline(1.5,.5)(1.5,-.5) \qline(1.5,-.5)(.5,-1) \qline(0,0)(0.5,-1)
\rput(1.75,.5){$\scriptstyle{3}$} \rput(-.5,-1){$H_1$}
}
\rput(3.5,1.5){%
\psdots*[](0,0)(.5,1)(1.5,1)(2,0)(1.5,-1)(.5,-1)(-1,0)
\qline(0,0)(.5,1) \qline(.5,1)(1.5,1) \qline(1.5,1)(2,0) \qline(2,0)(1.5,-1) \qline(1.5,-1)(.5,-1) \qline(0,0)(.5,-1) \qline(-1,0)(0,0) 
\rput(-.5,-1){$H_2$} \rput(1.5,1.25){$\scriptstyle{3}$}
}
\rput(7.5,1.5){%
\psdots*[](0,0)(.5,1)(1.5,1)(2,0)(1.5,-1)(.5,-1)(-1,0)
\qline(0,0)(.5,1) \qline(.5,1)(1.5,1) \qline(1.5,1)(2,0) \qline(2,0)(1.5,-1) \qline(1.5,-1)(.5,-1) \qline(0,0)(.5,-1) \qline(-1,0)(0,0)
\rput(-.5,-1){$H_3$} \rput(1.75,0){$\scriptstyle{3}$}
}
\rput(11.5,1.5){%
\psdots*[](0,0)(.5,1)(1.5,1)(2,0)(1.5,-1)(.5,-1)(-1,0)
\qline(0,0)(.5,1) \qline(.5,1)(1.5,1) \qline(1.5,1)(2,0) \qline(2,0)(1.5,-1) \qline(1.5,-1)(.5,-1) \qline(0,0)(.5,-1) \qline(-1,0)(0,0)
\rput(-.5,-1){$H_4$} \rput(-1,0.25){$\scriptstyle{3}$}
}
\endpspicture

Picture $4$ (continuation).\\
\end{center}

\begin{center}
\pspicture*[](-1,0)(14,8)
\rput(0,5.5){%
\psdots*[](0,0)(.5,0)(-.5,0)(0,.5)(0,-.5)(.5,.5)
\qline(-.5,0)(.5,0) \qline(0,-.5)(0,.5) \qline(0,0)(.5,.5)
\rput(.5,-.5){$H_5$}
}
\rput(3,5.5){%
\psdots*[](0,0)(.5,0)(-.5,0)(0,.5)(0,-.5)(0,1)
\qline(-.5,0)(.5,0) \qline(0,-.5)(0,1)
\rput(0,-1){$H_6\;\; (b\geq 2)$} \rput(-.25,1){${\scriptstyle b}$}
}
\rput(5,5.5){%
\psdots*[](0,0)(0,1)(0,-1)(1,0)(2,0)(4,0)(4,1)(4,-1)(4.5,1)
\qline(0,1)(0,-1) \qline(4,1)(4,-1) \qline(0,0)(2.9,0) \qline(3.1,0)(4,0) \qline(4.5,1)(4,1)
\psline[linestyle=dotted]{}(3,.2)(3,-.2)
\rput(2,-2){$H_7\;\; (0\leq n \leq 3)$} \rput(2,-.5){$\underbrace{\qquad\qquad\qquad\quad}_{{\displaystyle n \; }\mbox{\small{vertices}}}$}
}
\rput(11,5.5){%
\psdots*[](0,0)(.5,0)(.5,.5)(.5,-.5)(0,-1)(0,.5)(0,-.5)
\qline(0,-1)(0,.5) \qline(.5,-.5)(.5,.5) \qline(0,0)(.5,0)
\rput(-.25,-1){$\scriptstyle{3}$} \rput(1,-1){$H_8$}
}
\rput(1,1.5){%
\psdots*[](0,0)(.5,0)(1,0)(-.5,0)(-1,0)(-1.5,0)(0,-.5)(0,-1)
\qline(0,0)(0,-1) \qline(-1.5,0)(1,0)
\rput(1,-1){$H_9$}
}
\rput(5.5,1.5){%
\psdots*[](0,0)(.5,0)(1,0)(1.5,0)(2,0)(-.5,0)(-1,0)(-1.5,0)(0,-.5)
\qline(0,0)(0,-.5) \qline(-1.5,0)(2,0)
\rput(1,-1){$H_{10}$}
}
\rput(9.5,1.5){%
\psdots*[](0,0)(.5,0)(1,0)(1.5,0)(2,0)(2.5,0)(3,0)(-.5,0)(-1,0)(0,-.5)
\qline(0,0)(0,-.5) \qline(-1,0)(3,0)
\rput(3,.25){$\scriptstyle{b}$} \rput(2,-1){$H_{11}\;\; (b\geq 2)$}
}
\endpspicture

Picture $4$ (continuation).\\
\end{center}

{\it Proof:} The argument is essentially the same as in \cite{NikulinII}, Theorem~1.1.8.\\

Let $\Gamma$ be the associated graph of a minimal Lanner log canonical system with at least $4$ elements. If $\Gamma$ contains a double edge, then $\Gamma$ has exactly $4$ elements. It is obtained from $Q_8$, $Q_9$, $Q_{10}$ or $Q_{11}$ by adding a vertex. One checks that $\Gamma$ is one of $H_{8,1}$, ..., $H_{11,6}$ (Picture~4).\\

Assume that $\Gamma$ has no multiple edges. Since $\Gamma$ can not be a cycle, it is a tree, unless it is obtained from the cycle $\Gamma_5$ by adding one or two more vertices. If $\Gamma$ is a cycle with two more vertices, then $\Gamma$ is obtained from $Q_7$, $Q_{12}$, $Q_{13}$, $Q_{14}$ by adding an extra vertex. So, we get $H_{0,7,1}$, ..., $H_{0,14,1}$ (Picture~4).\\

If $\Gamma$ is a cycle with one more vertex, then $\Gamma$ is one of $H_{0,1}$, $H_{0,2}$, $H_{0,3}$, $H_{0,4}$ or one of $H_0$, $H_1$,..., $H_4$ (Picture~4).\\

Hence we may assume that $\Gamma$ is a tree without multiple edges. So, it should be one of $H_5$, ... $H_{11}$.\\

Note that one can not obtain a Lanner log canonical system by blowing up from Lanner systems with associated graphs $H_{8,1}$, ..., $H_{0,14,1}$, $H_{0,4}$, $H_{1}$, ..., $H_4$. {\it QED}\\

{\bf Corollary A.1.9 (cf. \cite{NikulinII}, Corollary 1.1.9).} {\it If $\Gamma$ is the associated graph of a Lanner log canonical system with at least $11$ elements, then $\Gamma$ is either a tree, or a cycle, or a cycle with one extra vertex and an edge. It has no multiple edges and every vertex of $\Gamma$ with weight $1$ has valency at most $2$. If $\Gamma$ is a cycle, then the Lanner system is obtained by a sequence of blow-ups of pairs of elements from $G_1(1,1;2)$, $G_1(1,2;2)$, $G_1(1,3;2)$, $G_3(1,1,b)$, where $b\geq 2$, or $G_6(1,1,1)$.}\\

{\it Proof:} This follows from Corollary A.1.5 and Corollary A.1.7 as in the proof of \cite{NikulinII}, Corollary 1.1.9. {\it QED}\\

{\bf Lemma A.3.0 (\cite{NikulinI}, \cite{AlexeevNikulin}).} {\it Let $X$ be a normal projective surface with log canonical singularities such that $-K_X$ is ample. Let $f\colon Y \to X$ be the minimal resolution of singularities. Then the set of all integral curves on $Y$ with negative self-intersection is finite. If $X$ is rational, then this set forms an at most hyperbolic system of vectors.\\

If $X$ is rational and different from ${\mathbb P}^2$, ${\mathbb P}^1\times {\mathbb P}^1$, the blow-up of ${\mathbb P}^2$ at one point and the contraction of the section with negative self-intersection on a Hirzebruch surface ${\mathbb F}_n={\mathbb P}(\mathcal O \oplus \mathcal O (n))\to {\mathbb P}^1$, $n\geq 2$, then this system is hyperbolic.}\\

{\it Proof:} The finiteness of the number of integral curves on $Y$ with negative self-intersection is proven in \cite{NikulinI}, Lemma 2.1 (i) (the proof goes through literally in our case). See also \cite{AlexeevNikulin}. Note that for an integral curve $C$ on $Y$ with negative self-intersection, either $C$ is contracted by $f$ or $C$ is a smooth rational $(-1)$-curve.\\

Let $\rho (Y)$ denote the Picard number of $Y$.\\

Recall (see \cite{KollarKovacs}, Corollary 2.1.4) that the Mori cone of $Y$ is generated by integral curves on $Y$ with negative self-intersection, unless $\rho (Y)=1$ (i.e. $Y=X={\mathbb P}^2$) or $\rho (Y)=2$ (i.e. either $Y=X$ is the blow-up of ${\mathbb P}^2$ at one point, or $Y=X={\mathbb P}^1\times {\mathbb P}^1$, or $X$ is the contraction of the section with negative self-intersection on a Hirzebruch surface $Y={\mathbb F}_n={\mathbb P}(\mathcal O \oplus \mathcal O (n))\to {\mathbb P}^1$, $n\geq 2$, or $X$ is the contraction of the section with negative self-intersection of a fibration $Y={\mathbb P}_E(\mathcal E)\to E$ for a normalized rank $2$ locally free sheaf $\mathcal E$ on a smooth elliptic curve $E$ with the parameter $e=-deg(\mathcal E)\geq 2$).\\

Hence the lattice determined by the intersection matrix of all integral curves on $Y$ with negative self-intersection has signature $(1,\rho (Y) - 1)$ by the Hodge Index theorem, unless 
\begin{itemize}
\item $Y={\mathbb P}^2$ or $Y={\mathbb P}^1\times {\mathbb P}^1$ (and there are no integral curves on $Y$ with negative self-intersection), or 
\item $Y={\mathbb F}_n$, $n\geq 1$, is a Hirzebruch surface (and there is a unique integral curve on $Y$ with negative self-intersection - the $(-n)$-section of the fibration ${\mathbb F}_n\to {\mathbb P}^1$), or
\item $X$ is not rational.
\end{itemize}

{\it QED}\\

Let $X$ be a rational normal projective surface with log canonical singularities such that $-K_X$ is ample. Let $f\colon Y \to X$ be the minimal resolution of singularities and $\mathcal V$ the at most hyperbolic system of vectors formed by the integral curves on $Y$ with negative self-intersection.\\

By \cite{Fujisawa}, Corollary~1.9, all singularities of $X$ are rational. Hence by \cite{KM}, Theorem~4.7 the associated weighted graph of the subsystem ${\mathcal V}_{\geq 2}$ of $\mathcal V$ formed by the integral curves on $Y$ with self-intersection at most $-2$ is a tree without multiple edges. ${\mathcal V}_{\geq 2}$ is an elliptic log canonical system (see \cite{NikulinII}, Lemma~2.1.2). Hence the connected components of its associated weighted graph are ${\Gamma}_1$, ${\Gamma}_2$, ${\Gamma}_3$ or ${\Gamma}_4$.\\

{\bf Proposition A.3.4 (cf. \cite{NikulinII}, Proposition 1.3.4).} {\it Let $X$ be a normal projective surface with log canonical singularities such that $-K_X$ is ample. Assume that $X$ is rational and different from ${\mathbb P}^2$, ${\mathbb P}^1\times {\mathbb P}^1$, the blow-up of ${\mathbb P}^2$ at one point and the contraction of the section with negative self-intersection on a Hirzebruch surface ${\mathbb F}_n={\mathbb P}(\mathcal O \oplus \mathcal O (n))\to {\mathbb P}^1$, $n\geq 2$.\\

Let $f\colon Y \to X$ be the minimal resolution of singularities and $\mathcal V$ the hyperbolic system of vectors formed by the integral curves on $Y$ with negative self-intersection.\\ 

If $\mathcal W\subset \mathcal V$ is a Lanner log canonical subsystem with at least $4$ elements, then $\mathcal W$ is obtained by a sequence of blow-ups from a Lanner log canonical system with the Gram matrix\\

$G_1(1,1;2)$, $G_1(1,2;2)$, $G_1(1,3;2)$, $G_1(1,1;3)$, $G_1(2,2;3)$, $G_1(2,2;4)$, $G_2(1,1,b)$, or\\

$G_3(1,1,b)$, $G_4(1,1,1)$, $G_6(1,1,1)$, where $b\geq 2$.}\\

{\it Proof:} We need to exclude the other Lanner log canonical systems listed in Theorem~A.1.8~(b). By the preceding remark $\mathcal V$ can not contain a subsystem with the associated graph $H_{8,1}$, ..., $H_{0,14,1}$, $H_{0,4}$, $H_1$, ..., $H_4$, as well as $H_{0,1}$, $H_{0,2}$, $H_{0,3}$.\\

For the remaining cases we follow the argument of Nikulin in \cite{NikulinII}, Proposition~1.3.4. It goes through and shows that $\mathcal W$ can not be obtained by a sequence of blow-ups from a Lanner log canonical system with the associated graph among $H_{0,1}$, $H_{0,2}$, $H_{0,3}$, $H_0$, $H_5$, ..., $H_{11}$ or from a Lanner log canonical system with the Gram matrix among\\

$G_2(1,4,b)$, $G_2(4,1,b)$, $G_3(1,4,b)$, $G_3(1,9,b)$, $G_4(1,1,4)$, $G_4(4,1,4)$, $G_5(4,1,1)$,\\

$G_5(4,1,4)$, $G_5(9,1,1)$, $G_6(1,1,4)$, $G_7(1,4,2)$, $G_7(4,1,4)$, $G_8(1,4,b)$, as well as among\\

$G_2(2,2,b)$, $G_3(2,2,b)$, $G_3(3,3,b)$, $G_5(2,2,2)$, $G_6(2,2,2)$, $G_6(3,3,3)$, $G_7(2,2,4)$,\\

$G_7(2,2,6)$, $G_8(2,2,b)$, $G_9(3,3,4)$, $G_9(2,2,b)$.\\

Note that when all proper subsystems of a Lanner log canonical system with the associated graph $H_{0,1}$, $H_{0,2}$ or $H_{0,3}$ are elliptic, one can not construct a different Lanner log canonical system out of it by blowing up. {\it QED}\\

The next Theorem is the main result of this Appendix. Its proof is literally the same as the proof of \cite{NikulinII},~Theorem~2.1.1.\\

{\bf Theorem A (cf. \cite{NikulinII}, Theorem~2.1.1).} {\it Let $X$ be a rational normal projective surface with log canonical singularities such that $-K_X$ is ample. Let $f\colon Y \to X$ be the minimal resolution of singularities and $\mathcal V$ the at most hyperbolic system of vectors formed by the integral curves on $Y$ with negative self-intersection. Then $\mathcal V$ is log canonical.}\\

\section*{Acknowledgement.}
We thank Beijing International Center for Mathematical Research and Peking University for providing excellent working conditions in which this note was prepared. We owe a lot to the works \cite{NikulinII} by Nikulin and \cite{Alexeev} by Alexeev. Paper \cite{Manoj} was also beneficial to us. We thank Manoj Verma for sending it to us.

\bibliographystyle{ams-plain}

\bibliography{FujinoConjectureForSurfaces}

\end{document}